\documentclass[a4paper]{article} 

\usepackage{a4wide}
\usepackage{authblk}
\usepackage{amsmath,amssymb,amsthm}
\usepackage{latexsym}
\usepackage{xcolor}
\usepackage{graphicx}
\usepackage{natbib}

\def\N{\mathbb{N}}

\def\R{\mathbb{R}}
\newcommand{\bracket}[1]{\left( #1 \right)}
\newcommand{\tri}{\triangle}
\newcommand{\PSR}[1]{#1_{\mathrm{PS}}}
\newcommand{\PST}{\PSR{\tri}}
\newcommand{\CTR}[1]{#1_{\mathrm{CT}}}
\newcommand{\CTT}{\CTR{\tri}}
\newcommand{\splspace}{\mathbb{S}}
\newcommand{\polspace}{\mathbb{P}}
\newcommand{\bloss}{\mathfrak{B}}
\newcommand{\Tsym}{T_{\mathrm{sym}}}
\newcommand{\nsym}{n_{\mathrm{sym}}}
\newcommand{\map}[1]{\mathcal{#1}}

\newtheorem{theorem}{Theorem}
\newtheorem{example}[theorem]{Example}
\newtheorem{remark}[theorem]{Remark}

\providecommand{\keywords}[1]{\noindent\small \textit{Keywords:} #1}

\begin{document}

%%%%%%%%%%%%%%%%%%%%%%%%%%%%%%%%%%%%%%%%%%%%%%%%%%%%%%%%%%%%%%%%%%%%%

\title{Extraction and application of super-smooth cubic B-splines over triangulations}

\author[$\dagger$,$\ddag$]{Jan Gro\v{s}elj}
\affil[$\dagger$]{Faculty of Mathematics and Physics, University of Ljubljana, Slovenia}
\affil[$\ddag$]{Institute of Mathematics, Physics and Mechanics, Ljubljana, Slovenia}
%\ead{jan.groselj@fmf.uni-lj.si}
\author[$\star$]{Hendrik Speleers}
\affil[$\star$]{Department of Mathematics, University of Rome Tor Vergata, Italy}
%\ead{speleers@mat.uniroma2.it}

\date{}

\maketitle

\begin{abstract}\noindent
The space of $C^1$ cubic Clough--Tocher splines is a classical finite element approximation space over triangulations for solving partial differential equations. However, for such a space there is no B-spline basis available, which is a preferred choice in computer aided geometric design and isogeometric analysis. A B-spline basis is a locally supported basis that forms a convex partition of unity. In this paper, we explore several alternative $C^1$ cubic spline spaces over triangulations equipped with a B-spline basis. They are defined over a Powell--Sabin refined triangulation and present different types of $C^2$ super-smoothness. The super-smooth B-splines are obtained through an extraction process, i.e., they are expressed in terms of less smooth basis functions. These alternative spline spaces maintain the same optimal approximation power as Clough--Tocher splines. This is illustrated with a selection of numerical examples in the context of least squares approximation and finite element approximation for second and fourth order boundary value problems.

\medskip
\end{abstract}

%%%%%%%%%%%%%%%%%%%%%%%%%%%%%%%%%%%%%%%%%%%%%%%%%%%%%%%%%%%%%%%%%%%%%

\section{Introduction}

Splines over triangulations are a common tool for bivariate approximation. They are inextricably linked to the finite element analysis, where they appear as triangular elements and are used in the process of finding an approximate solution to a partial differential equation. For many standard problems, it is sufficient to work with continuous (or even discontinuous) splines, which are easy to represent. Higher order problems, however, require smoother elements, a prominent example of which is the $C^1$ cubic Clough--Tocher element \citep{ct_clough_65}.

The advent of isogeometric analysis has put a particular focus on the representation of approximations in terms of spline basis functions. The driving idea of this method is to use the same spline model for design and analysis \citep{iga_cottrell_09}. In computer aided geometric design, B-spline functions are a standard tool for representing geometric objects, and here we understand these functions in a general sense, as functions that are locally supported and form a convex partition of unity. 
Furthermore, B-splines commonly have a high order of smoothness, which has proven to be beneficial in achieving an excellent accuracy with respect to the number of degrees of freedom and in obtaining a better spectral behavior in comparison with less smooth finite elements; see, e.g., \cite{hughes_08,evans_09,bressan_19,sande_19,manni_22}.

In the context of isogeometric analysis, the potential of smooth spline spaces over triangulations has been tested in several papers \citep{iga_speleers_12,iga_jaxon_14,beirao_15,iga_speleers_15,liu_18,iga_zareh_19}, but is not fully exploited yet. Compared with tensor-product spline spaces, using triangulations for tiling the domain offers more flexibility in describing the geometry and performing local refinements, however, this comes at the expense of more complicated spline constructions. Over a general unstructured triangulation, the dimension of $C^1$ spline spaces depends on the triangulation's geometry \citep{lai_07}, and so even identifying the right number of basis functions may be a tedious endeavor, possibly subjected to numerical instabilities. For this reason, it is customary to refine triangulations with particular splits (e.g., the Clough--Tocher split of each triangle into three smaller ones) which provide some structure to the spline partition and allow the association of degrees of freedom with topological properties, namely the vertices, edges, and triangles of the triangulation.

The next step in an attempt to bring triangular elements closer to the ideas of isogeometric analysis is the establishment of suitable basis functions. Despite recent developments in B-spline techniques for finite elements \citep{ps_dierckx_97,ps5_speleers_10,ct3_speleers_10,ps_speleers_13,arg_groselj_21,ct_groselj_22}, the standard $C^1$ cubic Clough--Tocher element still lacks a B-spline representation, and it is not even clear whether such a representation exists. In the absence of it, we explore in this paper alternative $C^1$ cubic spline spaces over triangulations that are readily equipped with B-spline basis functions \citep{ps3_groselj_17}. They are defined over a Powell--Sabin refined triangulation, which may be seen as a continuation of the Clough--Tocher split. Moreover, we focus on $C^2$ super-smoothness properties and develop an extraction process that enables a simple transition from less to more smooth B-spline basis functions.

It has been observed by \cite{ps3_groselj_21} that the general $C^1$ cubic Powell--Sabin B-splines admit additional $C^2$ super-smoothness properties when applied to three-directional meshes and that further smoothness properties can be obtained by combining and reducing certain B-splines. The aim of this paper is to provide an extraction procedure that applies those techniques to unstructured triangulations wherever possible and to numerically demonstrate the benefits of the procedure in approximation methods. As the main framework, we consider B-spline representations over an initial triangulation, which can be, for approximation purposes, uniformly refined. In this way, we acquire a three-directional triangulation on each initial triangle, where the symmetries can be exploited to reduce the number of degrees of freedom and to impose local $C^2$ super-smoothness.

In general, obtaining maximal smoothness for splines over triangulations with respect to the polynomial degree has proven to be very difficult. The techniques that come the closest to achieving this goal are based on simplex splines. This generalization of univariate B-splines is currently undergoing an exciting development \citep[for some recent advances, see][]{neamtu_07,liu_07, schmitt_21,barucq_22,ws3_lyche_22} and has been considered for use in the finite element method and isogeometric analysis \citep{cao_19,wang_22}. However, their disadvantage appears to be a rather complicated underlying configuration of knot lines, which makes the assembly and (numerical) integration a cumbersome problem when setting up the linear systems.

In the recent literature, several spline constructions have also been proposed over unstructured quadrilateral partitions consisting of large structured parts. In that context it is a common practice to use standard (tensor-product) B-spline basis functions supported over the structured parts of such a partition and specifically tailored basis functions (with often lower smoothness) in the vicinity of the interfaces between the structured parts; see, e.g., \cite{karciauskas_16,toshniwal_17,wei_18,chan_18,kapl_19,weinmuller_22}. The latter basis functions are usually not without shortcomings, e.g., they might require locally higher degrees, might not be nonnegative, might not be $C^1$ smooth, and/or might negatively affect the approximation power of the resulting space.
In this paper, we reverse this approach by starting with basis functions that are designed for general unstructured partitions and have none of these shortcomings, and then we reduce the basis functions over the structured parts of the partition to achieve locally higher smoothness. This process is described by extraction matrices allowing a fast transition from one approximation space to the other.

More particularly, we begin with the full space spanned by the $C^1$ cubic Powell--Sabin B-splines introduced by \cite{ps3_groselj_17}; these are associated with the vertices and edges of the triangulation. Then, firstly, we reduce the B-splines associated with the edges, which results in the spline space introduced by \cite{ps3_speleers_15} and asymptotically (with respect to the refinement of the initial triangulation) diminishes the dimension of the full space by a factor $5/3$. Secondly, we combine the previously reduced edge B-splines over all triangles that possess a specific geometric symmetry. The resulting space is locally the same as the one introduced by \cite{ps3_groselj_21} and asymptotically diminishes the dimension of the full space by a factor $3$.

All three considered spaces reproduce cubic polynomials and have optimal approximation power, which we illustrate with a selection of numerical examples in the context of least squares approximation and finite element approximation for second and fourth order boundary value problems. Through these examples we also demonstrate that locally increased smoothness improves accuracy of the approximations with respect to the number of degrees of freedom, i.e., the reduced spaces prove to be computationally more efficient.

The remainder of the paper is organized as follows. In Section~\ref{sec:basis}, we present $C^1$ cubic Powell--Sabin B-splines and explain the extraction process. In Section~\ref{sec:dim}, we derive dimension formulas for the considered spaces and investigate their asymptotic behavior. In Section~\ref{sec:applications}, we describe the application of the extraction process to a selection of approximation methods and provide a few numerical examples. We conclude with some remarks in Section~\ref{sec:conclusion}.

\section{Basis functions}
\label{sec:basis}

Suppose $\tri$ is an arbitrary triangulation of a planar bounded domain $\Omega$, and $\tri^\ell$, $\ell \in \N$, a refinement of $\tri$ obtained by splitting each triangle of $\tri$ uniformly such that each edge of $\tri$ is split into $\ell$ intervals as shown in Figure~\ref{fig:triangulation}. In this section, we will consider different sets of $C^1$ cubic B-spline basis functions over the triangulation $\tri^\ell$.

\subsection{Powell--Sabin refinement}

\begin{figure}[t]
\centering
\includegraphics[width=\textwidth]{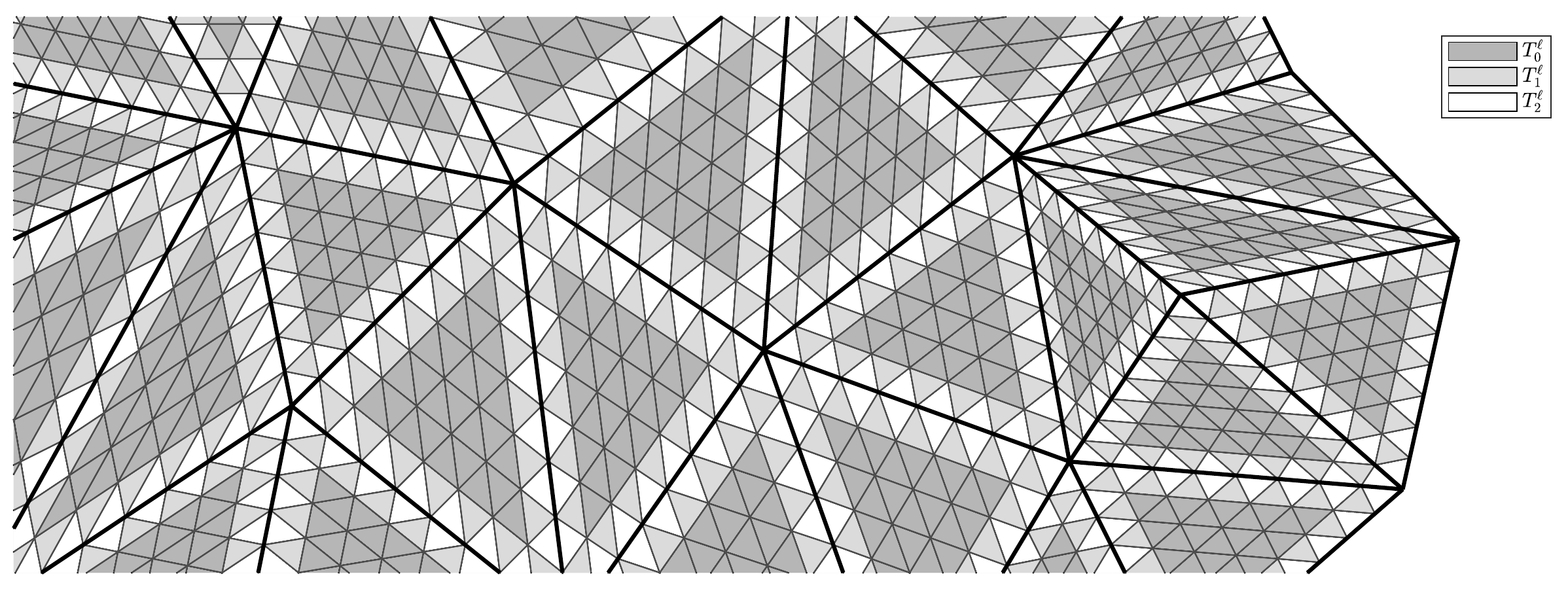}
\caption{A part of a triangulation $\tri$ (depicted by thicker lines) that is refined to the triangulation $\tri^\ell$ with $\ell = 8$. The color of a triangle in $\tri^\ell$ determines the set $T^\ell_r$, $r \in \{0, 1, 2\}$, to which the triangle belongs.}
\label{fig:triangulation}
\end{figure}

Let us represent the triangulation $\tri^\ell$ by the triplet $(V^\ell, E^\ell, T^\ell)$, where $V^\ell = \{ v_1, v_2, \ldots, v_{n_v^\ell} \}$ is the set of vertices, $E^\ell = \{ e_1, e_2, \ldots, e_{n_e^\ell} \}$ the set of edges, and $T^\ell = \{ t_1, t_2, \ldots, t_{n_t^\ell} \}$ the set of triangles. Note that $n_v^\ell$, $n_e^\ell$, and $n_t^\ell$ denote the number of elements in $V^\ell$, $E^\ell$, and $T^\ell$, respectively.
We further divide the set of triangles $T^\ell$ into three subsets $T^\ell_r$, $r = 0, 1, 2$. A triangle $t_k \in T^\ell$ belongs to $T^\ell_r$ if $r$ is the maximum number for which there exists an edge of $\tri$ that contains $r$ vertices of $t_k$; see Figure~\ref{fig:triangulation} for an illustration.

To introduce basis functions over $\tri^\ell$, we refine $\tri^\ell$ by applying the Powell--Sabin 6-split \citep{ps_powell_77}. For each triangle $t_k \in T^\ell$, we choose a split point $v_k^t$ inside $t_k$ such that for any neighboring triangle $t_{k'} \in T^\ell$ the line between $v_k^t$ and $v_{k'}^t$ intersects the interior of the edge $e_j \in E^\ell$ that is common to $t_k$ and $t_{k'}$. The intersection point is denoted by $v_j^e$. If $e_j$ is a boundary edge, the split point $v_j^e$ is the midpoint of the edge.

In our setting, a preferred position of $v_k^t$ is the barycenter of $t_k$, which is always a valid choice when $t_k \in T^\ell_0 \cup T^\ell_1$ as any two neighboring triangles in a uniform subtriangulation form a strictly convex quadrilateral. In general, however, choosing the barycenter of $t_k \in T^\ell_2$ could result in a degenerate or ill-poised refinement setting. On the other hand, choosing $v_k^t$ to be the incenter of $t_k$ guarantees the intersection property.

\begin{figure}[t!]
\centering
\includegraphics[width=\textwidth]{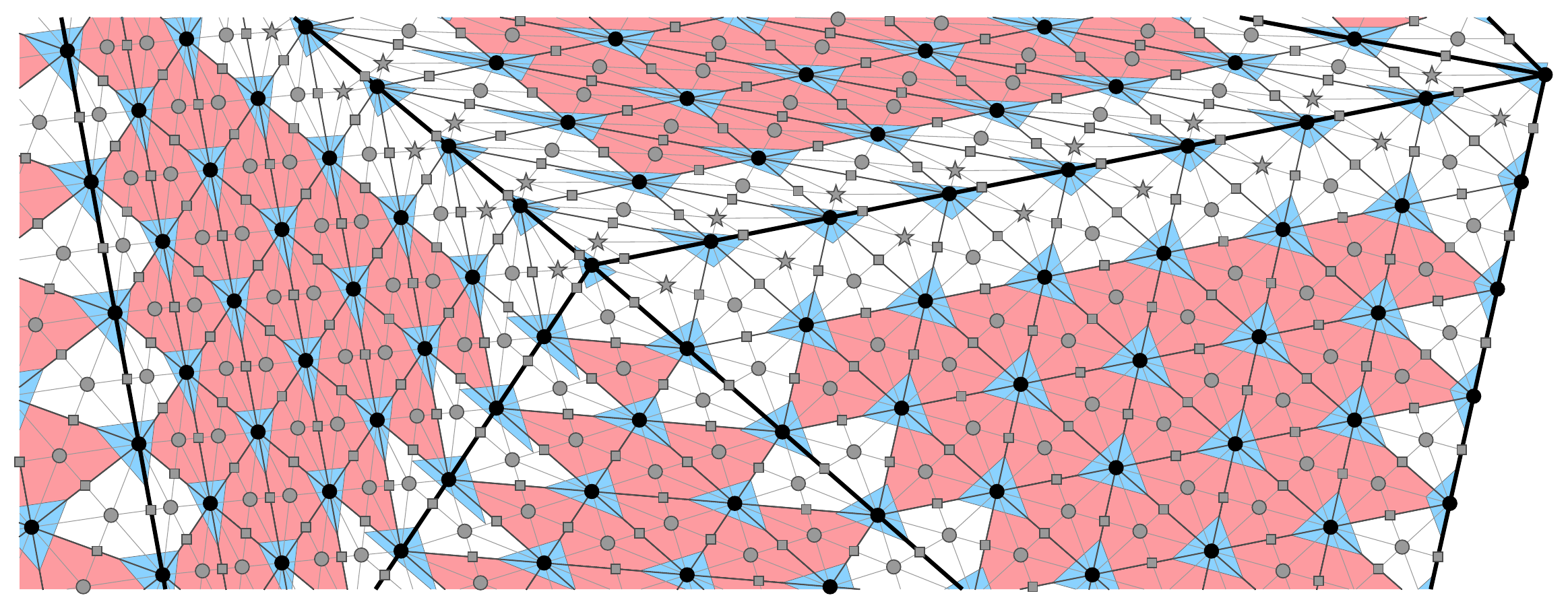}
\caption{An excerpt from the triangulation $\tri^\ell$ shown in Figure~\ref{fig:triangulation} and its Powell--Sabin refinement $\PST^\ell$. Each triangle of $\tri^\ell$ (depicted by dark gray lines) is split into six smaller triangles in $\PST ^\ell$ (depicted by light gray lines). The split points are colored gray. The ones inside triangles of $\tri^\ell$ are either barycenters or incenters, which are represented by circles or stars, respectively. The split points on the edges are marked by squares. The triangles that belong to $\Tsym^\ell$ are red shaded. The vertices of $\tri^\ell$ are represented by black points. For each $v_i \in V^\ell$, there is a triangle $q_i^v$ depicted in blue.}
\label{fig:pstriangulation}
\end{figure}

Each $t_k$ is split into six smaller triangles by connecting $v_k^t$ to the vertices of $t_k$ and to the split points on the edges of $t_k$. We denote the induced refinement of $\tri^\ell$ by $\PST^\ell$. An example of such a refinement is shown in Figure~\ref{fig:pstriangulation}.

Suppose a triangle $t_k$ of $\tri^\ell$ lies inside a uniform subtriangulation, i.e., if $v_{i_1}$, $v_{i_2}$, $v_{i_3}$ are the vertices of $t_k$, then its three neighboring triangles in $\tri^\ell$ are determined by the vertices
\begin{equation*}
v_{i_1}, v_{i_2}, v_{i_1} + v_{i_2} - v_{i_3}, \qquad
v_{i_2}, v_{i_3}, v_{i_2} + v_{i_3} - v_{i_1}, \qquad
v_{i_3}, v_{i_1}, v_{i_3} + v_{i_1} - v_{i_2}.
\end{equation*}
Suppose also that the split points of these four triangles are barycenters. We then say that $t_k$ has symmetric configuration. The set of all triangles in $\tri^\ell$ with symmetric configuration is denoted by $\Tsym^\ell$; see Figure~\ref{fig:pstriangulation} for an illustration.

\subsection{Full space}
\label{sec:full}

Over the triangulation $\PST^\ell$ of the domain $\Omega$, we consider the space $\splspace_0(\PST^\ell)$ consisting of $C^1$ functions such that the restriction to any triangle of $\PST^\ell$ is a bivariate polynomial of total degree less than or equal to $3$. The dimension of this space is $3 n_v^\ell + 4 n_e^\ell$; we refer the reader to \cite{ps3_groselj_17} for more details.

The space $\splspace_0(\PST^\ell)$ can be represented by basis functions $B_1, B_2, \ldots, B_m$, $m = 3 n_v^\ell$, associated with the vertices of $\tri^\ell$ and basis functions $B_1^0, B_2^0, \ldots, B_{m_0}^0$, $m_0 = 4 n_e^\ell$, associated with the edges of $\tri^\ell$.
We follow the approach from \cite{ps3_groselj_17} for their construction.

To each vertex $v_i \in V^\ell$, we assign a triangle $q_i^v$ that contains $v_i$ and the points $\frac{2}{3} v_i + \frac{1}{3} v_j^e$ and $\frac{2}{3} v_i + \frac{1}{3} v_k^t$ for each edge $e_j \in E^\ell$ and each triangle $t_k \in T^\ell$ attached to $v_i$; an illustration of such triangles can be found in Figure~\ref{fig:pstriangulation}. The vertices of $q_i^v$ happen to be Greville points \citep[see][for details]{ps3_groselj_17}.
To describe the other parameters in the construction, let us consider the edge $e_j \in E^\ell$ with the endpoints $v_i, v_{i'} \in V^\ell$ and the split point $v_j^e$. For each triangle of $\PST$ determined by the vertices $v_i$, $v_j^e$, $v_k^t$, we can choose a parameter $\sigma_{i,j,k} > 0$, and for a possible boundary edge of $\PST$ determined by the endpoints $v_i$, $v_j^e$, we can choose a parameter $\sigma_{i,j} > 0$. These yield the points $\frac{2}{3} v_i + \frac{1}{3} v_k^t + \frac{2}{3} \sigma_{i,j,k}(v_{i'} - v_i)$ and $\frac{2}{3} v_i + \frac{1}{3} v_j^e + \frac{2}{3} \sigma_{i,j}(v_{i'} - v_i)$, which play again the role of Greville points. An admissible choice that is assumed in the remainder of the paper is $\sigma_{i,j,k} = \sigma_{i,j} = \frac{1}{2}$.

\begin{figure}[t!]
\centering
\includegraphics[width=0.5\textwidth]{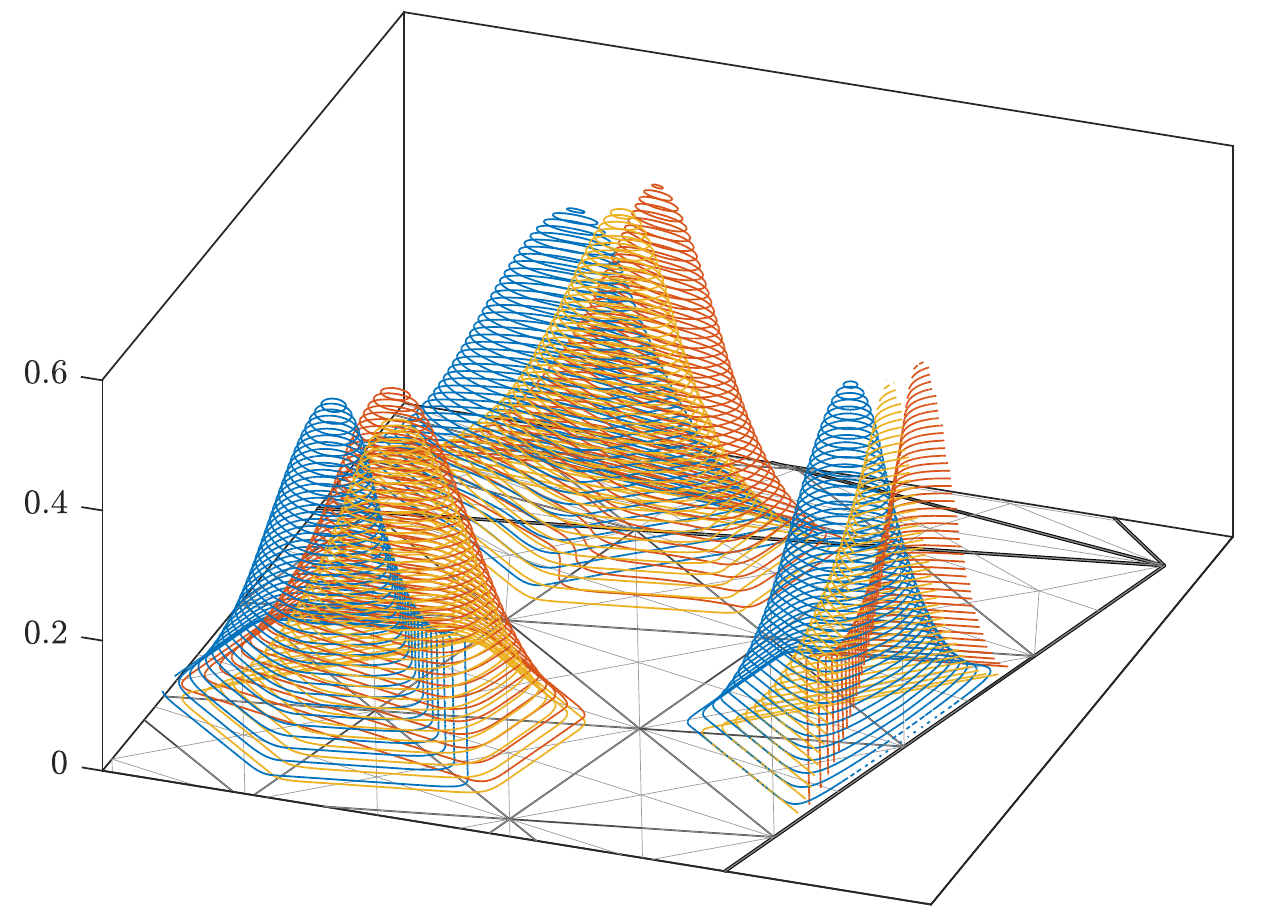}%
\includegraphics[width=0.5\textwidth]{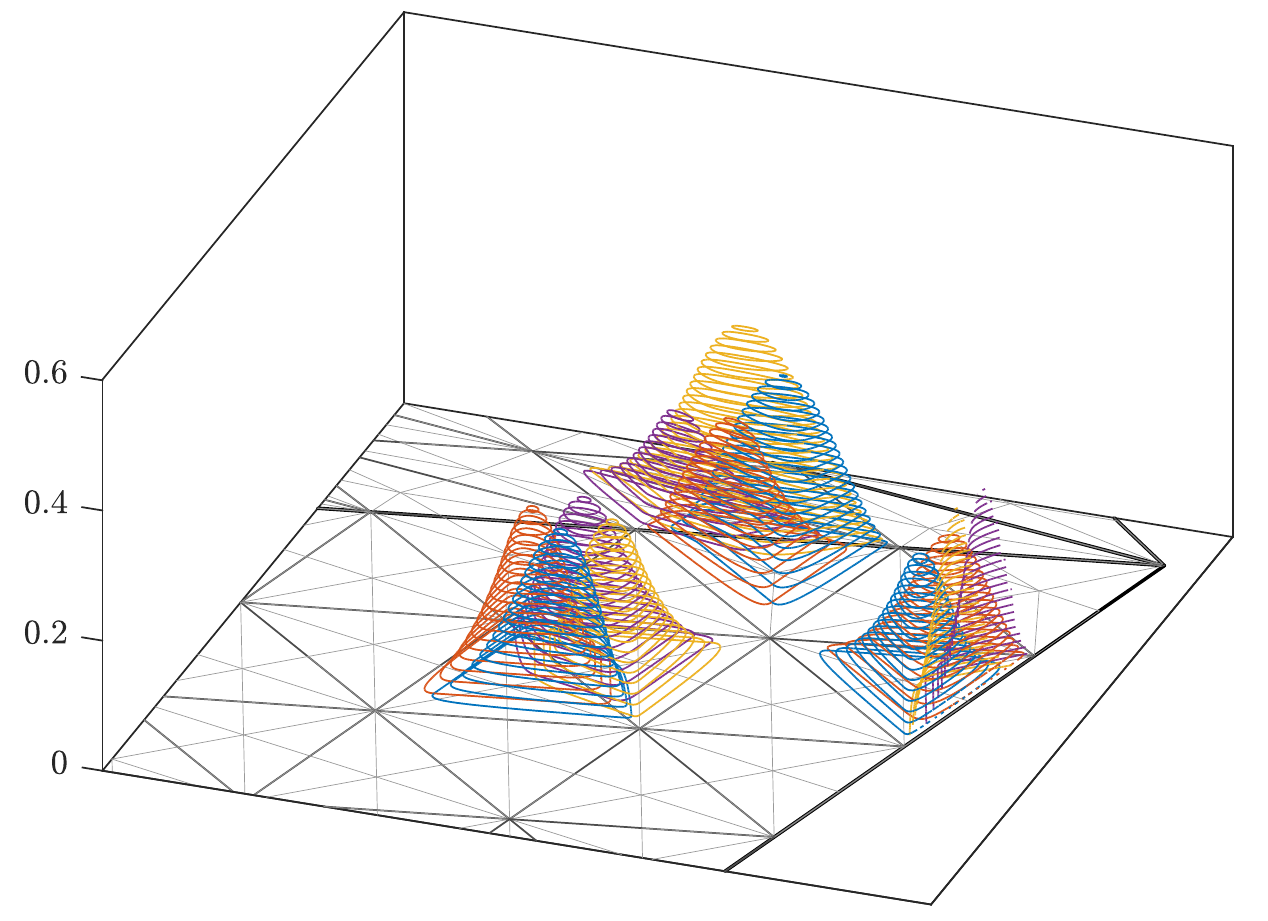}
\caption{Examples of basis functions $B_1, B_2, \ldots, B_m$ (left) and $B_1^0, B_2^0, \ldots, B_{m_0}^0$ (right) over $\PST^\ell$ shown in Figure~\ref{fig:pstriangulation}. The left picture shows the triplets of basis functions associated with three different types of vertices of $\tri^\ell$: a vertex encircled by triangles from $\Tsym^\ell$, a vertex on an inner edge of $\tri$, and a vertex on the boundary of the domain. The right picture shows the quartets of basis functions associated with three different types of edges of $\tri^\ell$: an edge shared by triangles from $\Tsym^\ell$, an edge lying on an inner edge of $\tri$, and an edge on the boundary of the domain.}
\label{fig:bsplines_0}
\end{figure}

Any triangle $q_i^v$ associated with $v_i \in V^\ell$ specifies three basis functions. 
Examples of such triplets of functions are shown in Figure~\ref{fig:bsplines_0} (left). These basis functions $B_1, B_2, \ldots, B_m$ possess the following super-smoothness properties; see \cite{ps3_groselj_17,ps3_groselj_21}.

\begin{theorem} \label{thm:smoothness_basis0}
The functions $B_1, B_2, \ldots, B_m$ possess the following super-smoothness properties.
\begin{enumerate}
\item For each $t_k \in T^\ell$, the functions are $C^2$ smooth at the split point $v_k^t$.
\item For each $t_k \in T^\ell$ and each edge $e_j$ of $t_k$, the functions are $C^2$ smooth across the edge connecting $v_j^e$ and $v_k^t$.
\item For each $t_k \in \Tsym^\ell$, the functions are $C^2$ smooth over the interior of the triangle $t_k$.
\end{enumerate}
\end{theorem}

Additionally, there are four basis functions associated with each edge of $\tri^\ell$. Examples of them are shown in Figure~\ref{fig:bsplines_0} (right). To distinguish between these basis functions $B_1^0, B_2^0, \ldots, B_{m_0}^0$, we assign them to triangles and boundary edges of $\PST^\ell$. Let $\Xi_0$ be a subset of $V^\ell \times E^\ell \times (T^\ell \cup \{ \emptyset \})$ consisting of $m_0$ elements. For each edge $e_j$ of $\tri^\ell$, it contains four triplets $(v_i, e_j, t_k)$, $(v_{i'}, e_j, t_k)$, $(v_i, e_j, t_{k'})$, $(v_{i'}, e_j, t_{k'})$ that identify
\begin{itemize}
\item four triangles of $\PST^\ell$ attached to $e_j$ if $e_j$ is not a boundary edge, or
\item two triangles and two boundary edges of $\PST^\ell$ if $e_j$ is a boundary edge ($t_{k'} = \emptyset$).
\end{itemize}
We denote by $\xi_0$ the bijection from $\Xi_0$ to $\{ 1, 2, \ldots, m_0 \}$. 

Contrary to $B_1, B_2, \ldots, B_m$, the functions $B_1^0, B_2^0, \ldots, B_{m_0}^0$ in general do not possess any additional super-smoothness properties. However, they complement the basis functions associated with the vertices in the following sense; see \cite{ps3_groselj_17}.

\begin{theorem} \label{thm:properties_basis0}
The functions $B_1, B_2, \ldots, B_m$ and the functions $B_1^0, B_2^0, \ldots, B_{m_0}^0$ form a convex partition of unity and a locally supported basis of $\splspace_0(\PST^\ell)$.
\end{theorem}

In the next sections, we consider two space reductions with additional local super-smoothness. The reduced basis functions are obtained through an extraction process, i.e., they are expressed in terms of less smooth basis functions.

\subsection{First space reduction}
\label{sec:reduction1}

Let us consider a subset $\Xi_1$ of $E^\ell \times (T^\ell \cup \{ \emptyset \})$ with $m_1 = 2 n_e^\ell$ elements. For each edge $e_j \in E^\ell$ it contains two pairs $(e_j, t_k)$, $(e_j, t_{k'})$ where $t_k$ and $t_{k'}$ are triangles of $\tri^\ell$ attached to $e_j$. If $e_j$ is a boundary edge, then $t_{k'} = \emptyset$. Suppose $\xi_1$ is a bijection from $\Xi_1$ to $\{ 1, 2, \ldots, m_1 \}$.

\begin{figure}[t!]
\centering
\includegraphics[width=0.5\textwidth]{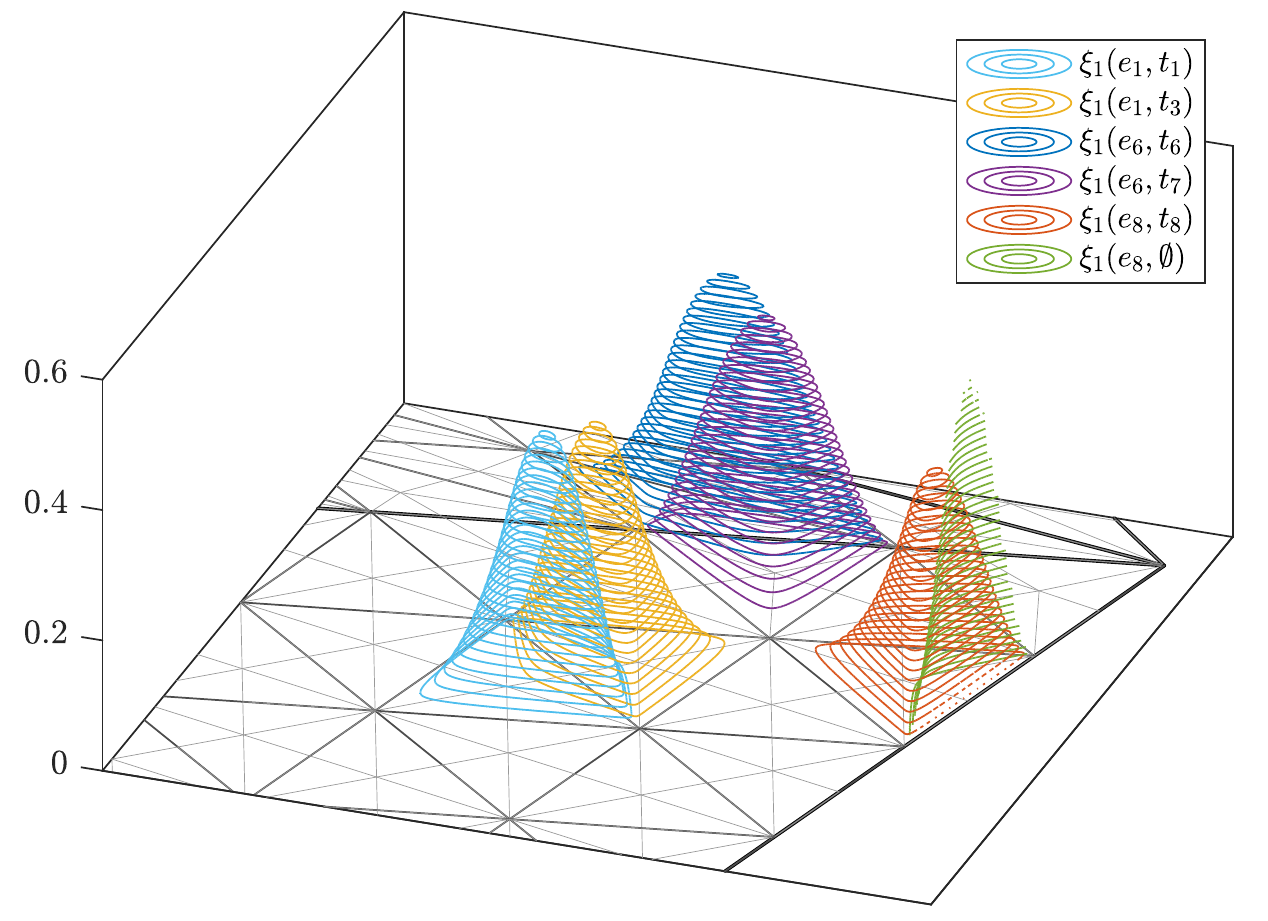}%
\includegraphics[width=0.5\textwidth]{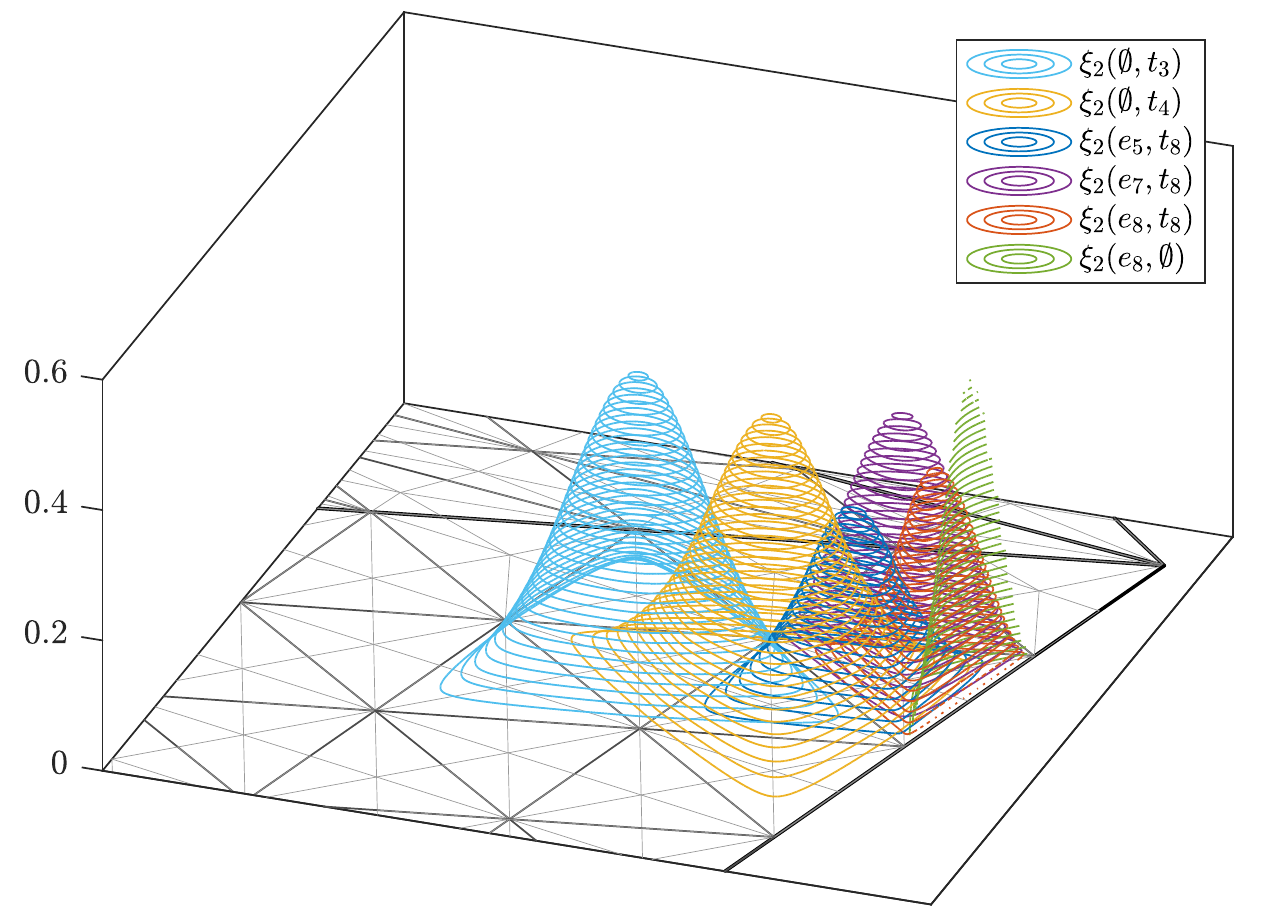}
\caption{Examples of basis functions $B_1^1, B_2^1, \ldots, B_{m_1}^1$ (left) and $B_1^2, B_2^2, \ldots, B_{m_2}^2$ (right) over $\PST^\ell$ shown in Figure~\ref{fig:pstriangulation}. The left picture shows the pairs of basis functions associated with three different types of edges of $\tri^\ell$: an edge shared by triangles from $\Tsym^\ell$, an edge lying on an inner edge of $\tri$, and an edge on the boundary of the domain. The right picture shows six functions: two basis functions spanning over four triangles of $\tri^\ell$ associated with two triangles from $\Tsym^\ell$, a basis function spanning over two triangles of $\tri^\ell$, one of which is in $\Tsym^\ell$ and the other is not, and three basis functions that remain unchanged after the second reduction (two of which are already shown on the left). For interpretation of the indices in the legends, see Figures \ref{fig:reduction1} and \ref{fig:reduction2}.}
\label{fig:bsplines12}
\end{figure}

Let $\map{B}_0: \Omega \rightarrow \mathbb{R}^{m_0}$ be a mapping that assigns to the point $p \in \Omega$ the column vector $[B_1^0(p), B_2^0(p), \ldots, B_{m_0}^0(p)]^T$. Furthermore, let $H_1 \in \mathbb{R}^{m_1 \times m_0}$ be the matrix with entries equal to $1$ at positions $(\xi_1(e_j, t_k), \xi_0(v_i, e_j, t_k))$ and $(\xi_1(e_j, t_k), \xi_0(v_{i'}, e_j, t_k))$ for each $(e_j, t_k) \in \Xi_1$, where $v_i, v_{i'} \in V^\ell$ are the endpoints of $e_j$. All other entries of $H_1$ are equal to $0$. The mapping $\map{B}_1: \Omega \rightarrow \mathbb{R}^{m_1}$ defined by 
\begin{equation*}
\map{B}_1(p) = H_1 \cdot \map{B}_0(p)
\end{equation*}
determines a new set of functions, which we denote by $B_1^1, B_2^1, \ldots, B_{m_1}^1$. Examples of them are shown in Figure~\ref{fig:bsplines12} (left). This procedure results in the reduced set of basis functions developed by \cite{ps3_speleers_15}. We refer the reader to \cite{ps3_groselj_17} for more details on the relation between the two sets of basis functions $\map{B}_0$ and $\map{B}_1$.

\begin{theorem} \label{thm:smoothness_basis1}
The functions $B_1^1, B_2^1, \ldots, B_{m_1}^1$ possess the following super-smoothness properties.
\begin{enumerate}
\item For each $t_k \in T^\ell$, the functions are $C^2$ smooth at the split point $v_k^t$.
\item For each $t_k \in T^\ell$ and each edge $e_j$ of $t_k$, the functions are $C^2$ smooth across the edge connecting $v_j^e$ and $v_k^t$.
\end{enumerate}
\end{theorem}

Let $\splspace_1(\PST^\ell)$ be the span of the functions $B_1, B_2, \ldots, B_m$ and the functions $B_1^1, B_2^1, \ldots, B_{m_1}^1$. The dimension of this space is $3 n_v^\ell + 2 n_e^\ell$.
From \cite{ps3_speleers_15} we also know the following properties.

\begin{theorem} \label{thm:properties_basis1}
The functions $B_1, B_2, \ldots, B_m$ and the functions $B_1^1, B_2^1, \ldots, B_{m_1}^1$ form a convex partition of unity and a locally supported basis of $\splspace_1(\PST^\ell)$.
\end{theorem}

\begin{figure}[t!]
\centering
\includegraphics[width=\textwidth]{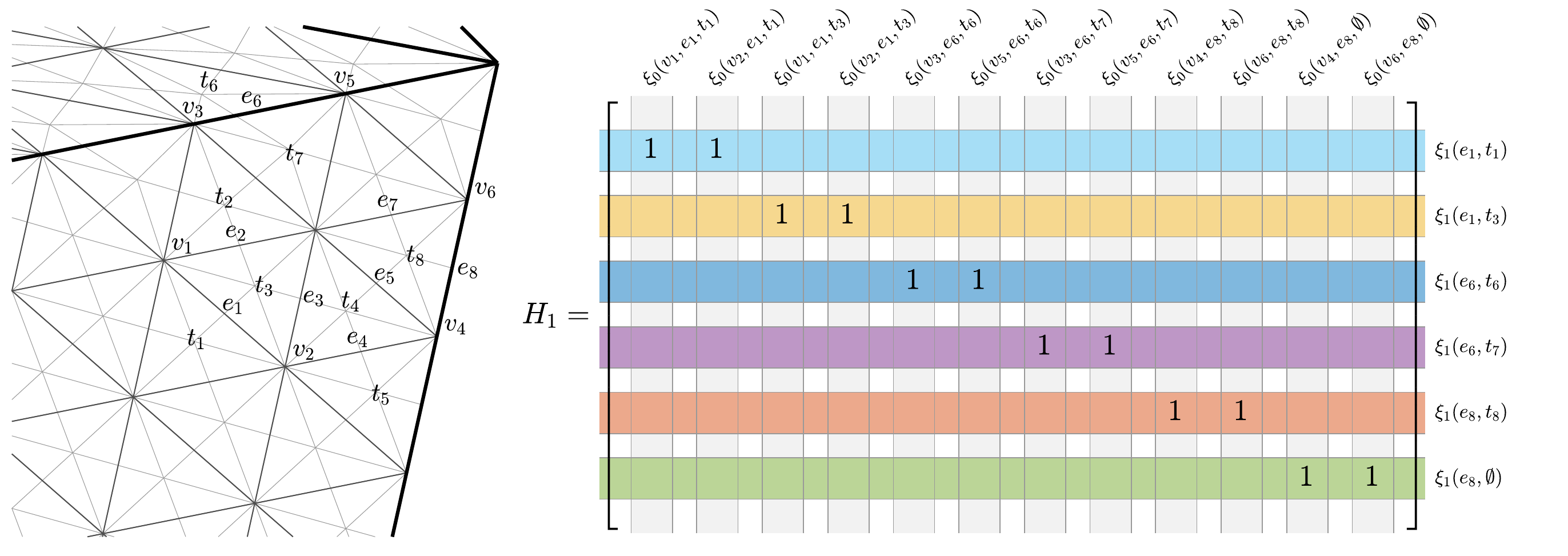}
\caption{An example of the first reduction. The basis functions shown in the left picture of Figure~\ref{fig:bsplines12} are obtained from the basis functions shown in the right picture of Figure~\ref{fig:bsplines_0}. An excerpt from the matrix $H_1$ contains the six rows that describe this reduction.}
\label{fig:reduction1}
\end{figure}

\begin{example}\label{ex:matrix_H1}
Figure~\ref{fig:reduction1} shows a few rows and columns of the matrix $H_1$ used for the reduction of the three quartets of basis functions depicted in Figure~\ref{fig:bsplines_0} (right) to the three pairs of basis functions depicted in Figure~\ref{fig:bsplines12} (left). Thanks to the properties of the original basis functions (see Theorem~\ref{thm:properties_basis0}), this matrix clearly confirms that the reduced basis functions form a convex partition of unity. Indeed, all its entries are nonnegative and the entries of each column sum to one. The local supports of the reduced basis functions follow from the sparsity of the matrix.
\end{example}

\begin{remark}
A matrix that describes a set of basis functions in terms of another set of (less smooth) basis functions is typically called an extraction matrix. This terminology was introduced by \cite{borden_11} and \cite{scott_11} to denote extraction of spline functions in terms of Bernstein polynomials. More recent examples of extraction matrices can be found in \cite{mdb_speleers19} and \cite{polar_speleers_21}.
\end{remark}

\subsection{Second space reduction}
\label{sec:reduction2}

Let $\Xi_2 \subset (E^\ell \cup \{ \emptyset \}) \times (T^\ell \cup \{ \emptyset \})$ be the union of $\{ \emptyset \} \times \Tsym^\ell$ and $\Xi_1 \cap (E^\ell \times (T^\ell \setminus \Tsym^\ell \cup \{ \emptyset \}))$, and let $m_2$ denote the number of elements of this set. Suppose $\xi_2$ is a bijection from $\Xi_2$ to $\{ 1, 2, \ldots, m_2 \}$.

Let $H_2 \in \mathbb{R}^{m_2 \times m_1}$ be the matrix with entries defined as follows. For each $t_k \in \Tsym^\ell$ and each edge $e_j$ of $t_k$, the entry at position $(\xi_2(\emptyset, t_k), \xi_1(e_j,t_k))$ is equal to $\frac{2}{3}$. For the second triangle $t_{k'} \in T^\ell$ attached to $e_j$, the entry at position $(\xi_2(\emptyset, t_k), \xi_1(e_j, t_{k'}))$ is equal to $\frac{1}{3}$. Moreover, if $t_{k'} \notin \Tsym^\ell$, the entry at position $(\xi_2(e_j, t_{k'}), \xi_1(e_j, t_{k'}))$ is equal to $\frac{2}{3}$ and the entry at position $(\xi_2(e_j, t_{k'}), \xi_1(e_j, t_k))$ is equal to $\frac{1}{3}$. Finally, for each pair $(e_j, t_k), (e_j, t_{k'}) \in \Xi_1$ such that both $t_k$ and $t_{k'}$ are in $T^\ell \setminus \Tsym^\ell \cup \{ \emptyset \}$, the entries at positions $(\xi_2(e_j,t_k), \xi_1(e_j, t_k))$ and $(\xi_2(e_j, t_{k'}), \xi_1(e_j, t_{k'}))$ are equal to $1$. All other entries of $H_2$ are equal to $0$. Note that if $\Tsym^\ell$ is empty, $H_2$ is a permutation of the identity matrix.

The mapping $\map{B}_2: \Omega \rightarrow \mathbb{R}^{m_2}$ defined by 
\begin{equation*}
\map{B}_2(p) = H_2 \cdot \map{B}_1(p) 
\end{equation*}
determines a new set of functions, which we denote by $B_1^2, B_2^2, \ldots, B_{m_2}^2$. Examples of them are shown in Figure~\ref{fig:bsplines12} (right). For all $t_k\in \Tsym^\ell$, the new functions $B_{\xi_2(\emptyset, t_k)}^2$ coincide with the basis functions described by \cite{ps3_groselj_21}.

\smallskip
\begin{theorem} \label{thm:smoothness_basis2}
The functions $B_1^2, B_2^2, \ldots, B_{m_2}^2$ possess the following super-smoothness properties.%
\begin{enumerate}
\item For each $t_k \in T^\ell$, the functions are $C^2$ smooth at the split point $v_k^t$.
\item For each $t_k \in T^\ell$ and each edge $e_j$ of $t_k$, the functions are $C^2$ smooth across the edge connecting $v_j^e$ and $v_k^t$.
\item For each $t_k \in \Tsym^\ell$, the functions are $C^2$ smooth over the interior of the triangle $t_k$.
\end{enumerate}
\end{theorem}

Let $\splspace_2(\PST^\ell)$ be the span of the functions $B_1, B_2, \ldots, B_m$ and the functions $B_1^2, B_2^2, \ldots, B_{m_2}^2$.
The following properties can be deduced by combining results from \cite{ps3_groselj_21} and \cite{ps3_speleers_15}.

\begin{theorem} \label{thm:properties_basis2}
The functions $B_1, B_2, \ldots, B_m$ and the functions $B_1^2, B_2^2, \ldots, B_{m_2}^2$ form a convex partition of unity and a locally supported basis of $\splspace_2(\PST^\ell)$.
\end{theorem}

\begin{figure}[t!]
\centering
\includegraphics[width=\textwidth]{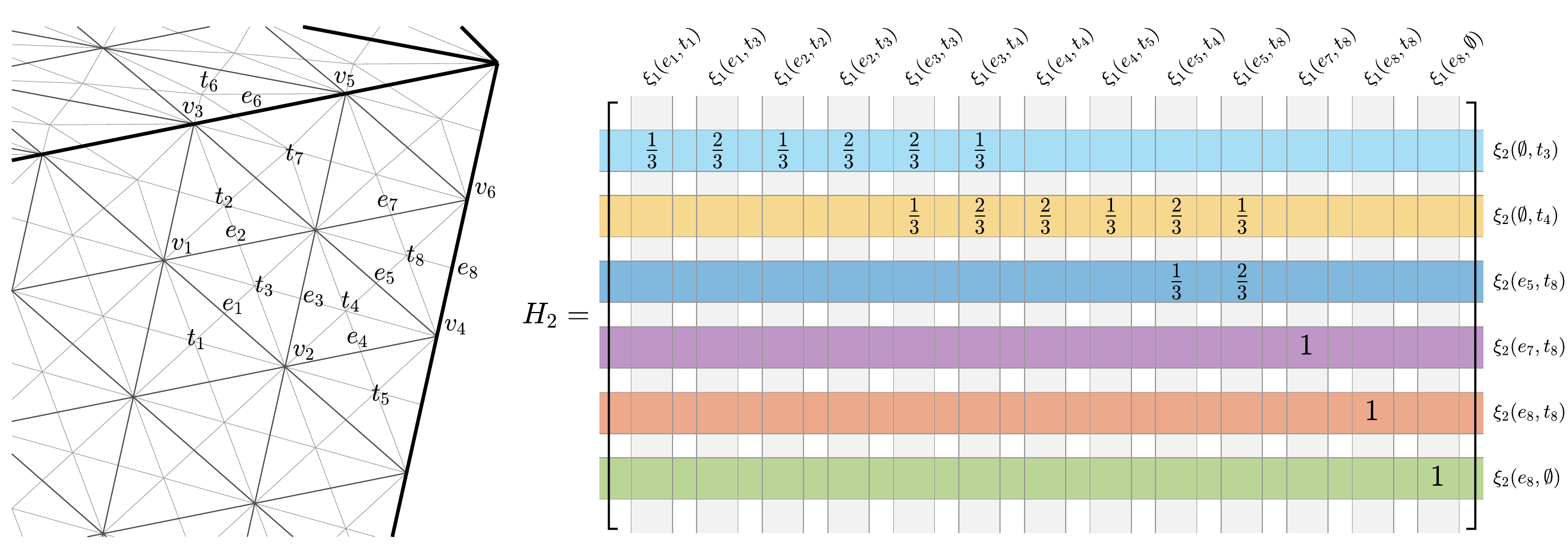}
\caption{An example of the second reduction. The basis functions shown in the right picture of Figure~\ref{fig:bsplines12} are produced from the basis functions shown in the left picture of Figure~\ref{fig:bsplines12} and a number of additional basis functions from the first reduction. An excerpt from the matrix $H_2$ contains the six rows that describe this reduction.}
\label{fig:reduction2}
\end{figure}

\begin{example} \label{ex:matrix_H2}
Figure~\ref{fig:reduction2} shows a few rows and columns of the matrix $H_2$ which determine the 6 reduced basis functions depicted in Figure~\ref{fig:bsplines12} (right). Thanks to the properties of the first reduced basis functions (see Theorem~\ref{thm:properties_basis1}), this matrix clearly confirms that the second reduced basis functions form a convex partition of unity. Indeed, all its entries are nonnegative and the entries of each column sum to one. The local supports of the second reduced basis functions follow from the sparsity of the matrix.
\end{example}

In Appendix~\ref{sec:polynomials} we show that the spline space $\splspace_2(\PST^\ell)$ contains all bivariate polynomials of total degree less than or equal to $3$. This implies that the space reduction presented here maintains the same optimal approximation power as the other spline spaces $\splspace_1(\PST^\ell) \subset \splspace_0(\PST^\ell)$.

\section{Dimension formulas}
\label{sec:dim}

The three spline spaces considered in the previous section differ in local super-smoothness, but all reproduce cubic polynomials and possess an optimal approximation order of four. In this section, we will analyze and compare their dimensions.

Let $\tri$ be an arbitrary triangulation consisting of $n_v$ vertices, $n_e$ edges, $n_{\partial e}$ boundary edges, and $n_t$ triangles.
Recall that $\tri^\ell$, $\ell \in \N$, is obtained from $\tri$ by splitting each triangle of $\tri$ uniformly such that each edge of $\tri$ is split into $\ell$ intervals. This refined triangulation has $n_v^\ell$ vertices, $n_e^\ell$ edges, $n_{\partial e}^\ell$ boundary edges, and $n_t^\ell$ triangles, where
\begin{align*}
n_v^\ell &= n_v + (\ell-1) n_e + \frac{1}{2}(\ell-2)(\ell-1) n_t,\\
n_e^\ell &= \ell\, n_e + \frac{3}{2}(\ell-1)\ell\, n_t,\\
n_{\partial e}^\ell &= \ell n_{\partial e} = 2\ell n_e - 3\ell n_t,\\
n_t^\ell &= \ell^2\, n_t.
\end{align*}
Moreover, let $\nsym^\ell$ be the number of triangles in $\Tsym^\ell$. For simplicity of exposition, we assume that the triangles in $\tri$ do not form locally uniform subtriangulations.

The dimension of the full spline space $\splspace_0(\PST^\ell)$ is given by
\begin{equation*}
\dim(\splspace_0(\PST^\ell)) = 3n_v^\ell + 4n_e^\ell 
= 3 n_v + (7\ell - 3) n_e + \frac{3}{2}(5\ell-2)(\ell-1) n_t.
\end{equation*}
Similarly, we deduce that the spline space after the first reduction $\splspace_1(\PST^\ell)$ has dimension
\begin{equation*}
\dim(\splspace_1(\PST^\ell)) = 3n_v^\ell + 2n_e^\ell 
= 3 n_v + (5\ell - 3) n_e + \frac{3}{2}(3\ell-2)(\ell-1) n_t.
\end{equation*}
Finally, the dimension of the second reduced spline space $\splspace_2(\PST^\ell)$ equals
\begin{align}
\dim(\splspace_2(\PST^\ell)) &= 3n_v^\ell \tag{corresponding to vertices} \\
   &\quad+ \nsym^\ell \tag{corresponding to symmetric triangles} \\
   &\quad + 3(n_t^\ell-\nsym^\ell) + n_{\partial e}^\ell, \tag{corresponding to asymmetric triangles}
\end{align}
where the number $\nsym^\ell$ depends on the geometry of $\PST^\ell$. In Example~\ref{ex:n_sym} we illustrate two extreme cases of $\nsym^\ell$. 

\begin{figure}[t!]
\centering
\includegraphics[width=0.5\textwidth]{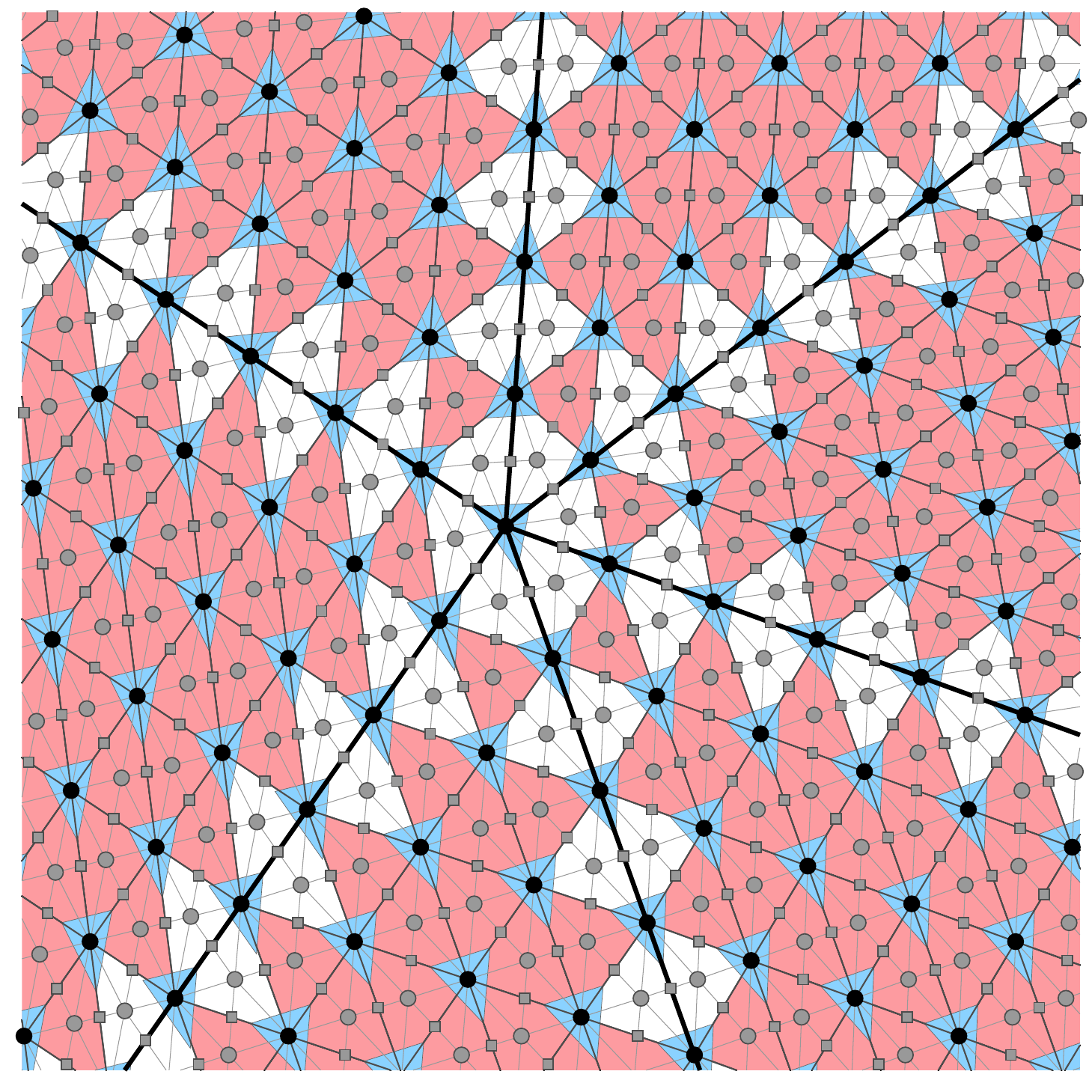}%
\includegraphics[width=0.5\textwidth]{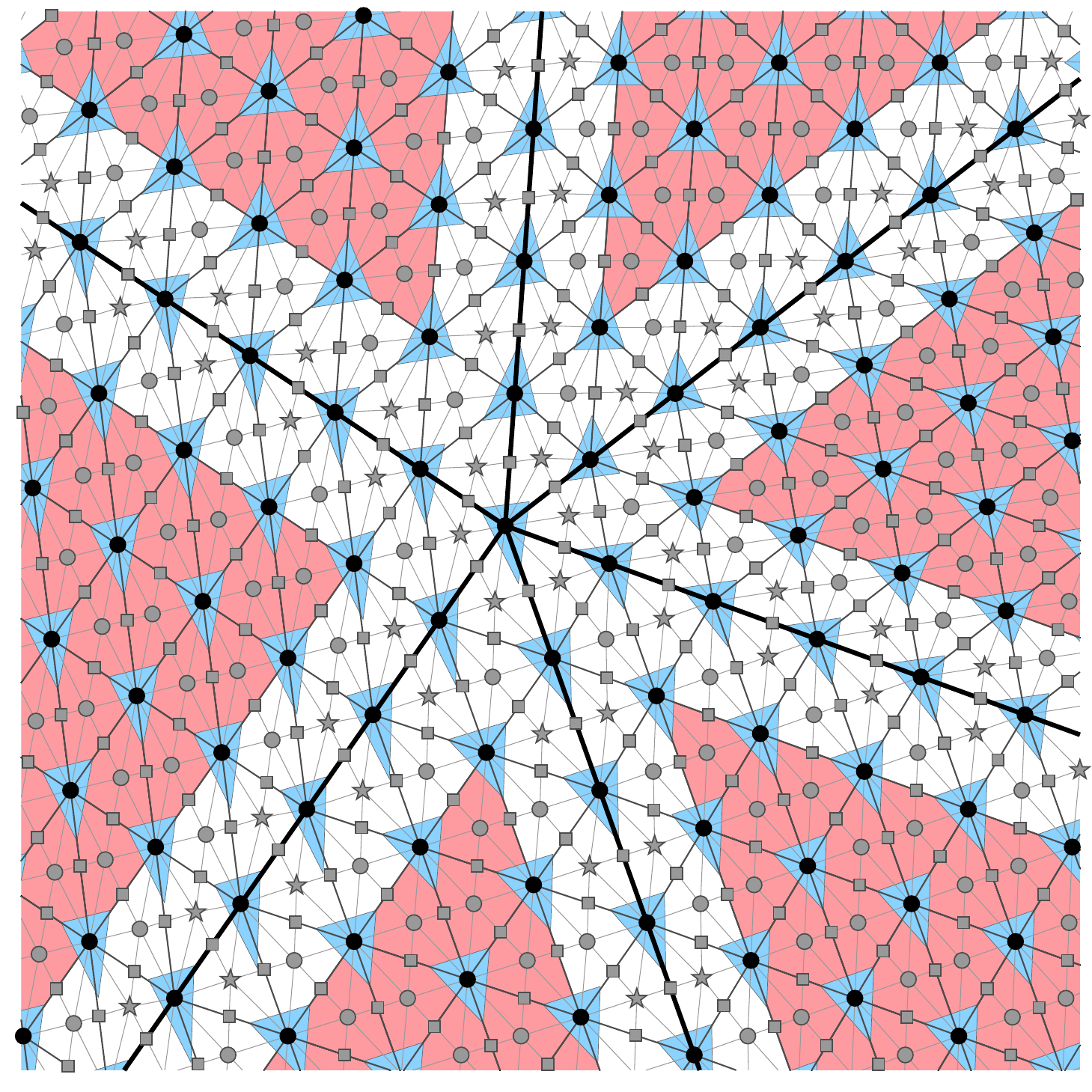}
\caption{Illustration of two extreme cases of $\Tsym^\ell$. The marks and color codes are the same as in Figure~\ref{fig:pstriangulation}.}
\label{fig:T_sym}
\end{figure}
\begin{example}\label{ex:n_sym}
If all the triangles in $\tri$ are refined as shown in the left picture of Figure~\ref{fig:T_sym}, then $\Tsym^\ell = T_0^\ell \cup T_1^\ell$ and the number of triangles in $\Tsym^\ell$ is maximal, namely
\begin{equation*}
\nsym^\ell = [\ell^2 - 3(\ell-1)] n_t, \quad \ell\geq2.
\end{equation*}
On the other hand, if all the triangles in $\tri$ are refined as shown in the right picture of Figure~\ref{fig:T_sym}, then $\Tsym^\ell = T_0^\ell$ and the number of triangles in $\Tsym^\ell$ is minimal, namely
\begin{equation*}
\nsym^\ell = (\ell-3)^2 n_t, \quad \ell\geq3.
\end{equation*}
In general, we have
\begin{equation*}
(\ell-3)^2 n_t \leq \nsym^\ell \leq [\ell^2 - 3(\ell-1)] n_t, \quad \ell\geq3.
\end{equation*}
An illustration of a general configuration can be found in Figure~\ref{fig:pstriangulation}.
\end{example}

In the following we fix $\ell\geq3$. We can rewrite the above dimension formula of $\splspace_2(\PST^\ell)$ as
\begin{equation*}
\dim(\splspace_2(\PST^\ell)) = 3n_v + (5\ell-3) n_e + \frac{3}{2}(3\ell^2-5\ell+2) n_t - 2\nsym^\ell,
\end{equation*}
and from the bounds for $\nsym^\ell$ provided in Example~\ref{ex:n_sym} we obtain
\begin{equation*}
 \dim(\splspace_2(\PST^\ell)) \geq 3n_v + (5\ell-3) n_e + \frac{1}{2}(5\ell^2-3\ell-6) n_t
\end{equation*}
and 
\begin{equation*}
 \dim(\splspace_2(\PST^\ell)) \leq 3n_v + (5\ell-3) n_e + \frac{1}{2}(5\ell^2+9\ell-30) n_t.
\end{equation*}

By comparing the different dimension formulas we deduce
\begin{align*}
 \dim(\splspace_0(\PST^\ell)) - \dim(\splspace_1(\PST^\ell)) &= 2\ell\, n_e + 3(\ell-1)\ell\, n_t  > 0,\\
 \dim(\splspace_0(\PST^\ell)) - \dim(\splspace_2(\PST^\ell)) &\geq 2\ell\, n_e + (5\ell^2-15\ell+18) n_t > 0.
\end{align*}
We also observe
\begin{align*}
 \lim_{\ell\rightarrow\infty} \frac{\dim(\splspace_0(\PST^\ell))}{\dim(\splspace_1(\PST^\ell))} &= \frac{15/2}{9/2} = \frac{5}{3},\\
 \lim_{\ell\rightarrow\infty} \frac{\dim(\splspace_0(\PST^\ell))}{\dim(\splspace_2(\PST^\ell))} &= \frac{15/2}{5/2} = 3.
\end{align*}
Thus, asymptotically speaking, the dimension drops by a factor $5/3$ and $3$, respectively, for the two reduced spaces without loss of approximation order.

\begin{remark}
The space of $C^1$ cubic Clough--Tocher splines is a well-known space over triangulations; see \cite{ct_clough_65,lai_07}. In this case, the triangulation $\tri^\ell$ is further refined to a triangulation, denoted by $\CTT^\ell$, where each triangle $t_k \in T^\ell$ is split into three smaller triangles by connecting a split point $v_k^t$ inside $t_k$ to the vertices of $t_k$. Let us denote the corresponding $C^1$ cubic spline space by $\splspace_0(\CTT^\ell)$. When the partitions $\CTT^\ell$ and $\PST^\ell$ are compatible, it has been shown by \cite{ps3_speleers_15} that $\splspace_0(\CTT^\ell)$ is a subspace of $\splspace_1(\PST^\ell)$ and, as a consequence, also of $\splspace_0(\PST^\ell)$. The dimension of $\splspace_0(\CTT^\ell)$ equals
\begin{equation*}
\dim(\splspace_0(\CTT^\ell)) = 3n_v^\ell + n_e^\ell 
= 3 n_v + (4\ell - 3) n_e + 3(\ell-1)^2 n_t.
\end{equation*}
It is asymptotically slightly larger than the dimension of $\splspace_2(\PST^\ell)$. More precisely,
\begin{equation*}
 \lim_{\ell\rightarrow\infty} \frac{\dim(\splspace_0(\CTT^\ell))}{\dim(\splspace_2(\PST^\ell))} = \frac{3}{5/2} = \frac{6}{5}.
\end{equation*}
Moreover, there is no B-spline basis available for the space $\splspace_0(\CTT^\ell)$. Partial B-spline results can be found in \cite{ct3_speleers_10} and \cite{ct3_lyche_18}.
\end{remark}

\section{Applications}
\label{sec:applications}

In this section, we provide numerical examples demonstrating the performance of the considered spline spaces in the context of least squares approximation and finite element approximation for second and fourth order boundary value problems. We start off by explaining how to take advantage of the reductions presented in Sections~\ref{sec:reduction1} and \ref{sec:reduction2}.

\subsection{Practical aspects of basis reduction}

Analogous to the mapping $\map{B}_0$ defined in Section~\ref{sec:full}, let us denote by $\map{B}: \Omega \rightarrow \R^m$ the mapping that assigns to the point $p \in \Omega$ the column vector $[B_1(p), B_2(p), \ldots, B_m(p)]^T$. Suppose $\map{S}_r: \Omega \rightarrow \R^{m+m_r}$, $r \in \{0, 1, 2\}$, is the mapping that assigns to $p \in \Omega$ the column vector obtained by combining the vectors $\map{B}(p)$ and $\map{B}_r(p)$. Then, $\map{S}_r$ determines a basis of $\splspace_r(\PST^\ell)$, and for each $s_r \in \splspace_r(\PST^\ell)$, there exists a column vector $c_r \in \R^{m+m_r}$ such that
\begin{equation}
\label{eq:spline}
s_r(p) = \map{S}_r(p)^T \cdot c_r.
\end{equation}
The mappings $\map{S}_r$ can be related by extended extraction matrices as
\begin{equation*}
\map{S}_1(p) = \widetilde{H}_1 \cdot \map{S}_0(p),
\quad
\map{S}_2(p) = \widetilde{H}_2 \cdot \map{S}_1(p),
\end{equation*}
where
\begin{equation*}
\widetilde{H}_1 =
\begin{bmatrix}
I & \\
& H_1
\end{bmatrix},
\quad
\widetilde{H}_2 =
\begin{bmatrix}
I & \\
& H_2
\end{bmatrix},
\end{equation*}
and $I \in \R^{m \times m}$ denotes the identity matrix.
Moreover, for $\alpha, \beta, \gamma \in \N_0$ such that $\alpha, \beta \leq \gamma$, let $\map{M}_r^{(\alpha,\beta,\gamma)}: \Omega \rightarrow \R^{(m+m_r) \times (m+m_r)}$ be the mappings defined by
\begin{equation}
\label{eq:splineprodmatrix}
\map{M}_r^{(\alpha,\beta,\gamma)}(p) = \bracket{D_x^{\alpha} D_y^{\gamma - \alpha} \map{S}_r(p)} \cdot \bracket{D_x^{\beta} D_y^{\gamma-\beta} \map{S}_r(p)}^T.
\end{equation}
Since
\begin{equation*}
\begin{aligned}
\map{M}_1^{(\alpha,\beta,\gamma)}(p) &= \widetilde{H}_1 \cdot \map{M}_0^{(\alpha,\beta,\gamma)}(p) \cdot \widetilde{H}_1^T, \\
\map{M}_2^{(\alpha,\beta,\gamma)}(p) &= \widetilde{H}_2 \cdot \map{M}_1^{(\alpha,\beta,\gamma)}(p) \cdot \widetilde{H}_2^T,
\end{aligned}
\end{equation*}
it is sufficient to compute $\map{M}_0^{(\alpha,\beta,\gamma)}(p)$. Then, $\map{M}_1^{(\alpha,\beta,\gamma)}(p)$ and $\map{M}_2^{(\alpha,\beta,\gamma)}(p)$ are easy to obtain by matrix multiplications. The considered mappings can be used to express matrices that typically arise in approximation methods. For example, as we demonstrate in the next sections, the mass and stiffness matrices can be acquired by elementwise integration of these mappings.

\subsection{Least squares approximation}

For a function $f \in L^2(\Omega)$, the least squares approximation $s_r$ in $\splspace_r(\PST^\ell)$ is given by the orthogonality condition
\begin{equation*}
\int_{\Omega} \map{S}_r(p) \cdot \bracket{f(p) - s_r(p)} \, \mathrm{d}p = 0,
\end{equation*}
where the integral is considered elementwise. By \eqref{eq:spline} and \eqref{eq:splineprodmatrix} this is the same as
\begin{equation*}
\bracket{ \int_\Omega \map{M}_r^{(0,0,0)}(p) \, \mathrm{d}p } \cdot c_r = \int_\Omega \map{S}_r(p) \cdot f(p) \, \mathrm{d}p.
\end{equation*}

\begin{figure}[t!]
\centering
\includegraphics[width=0.33\textwidth]{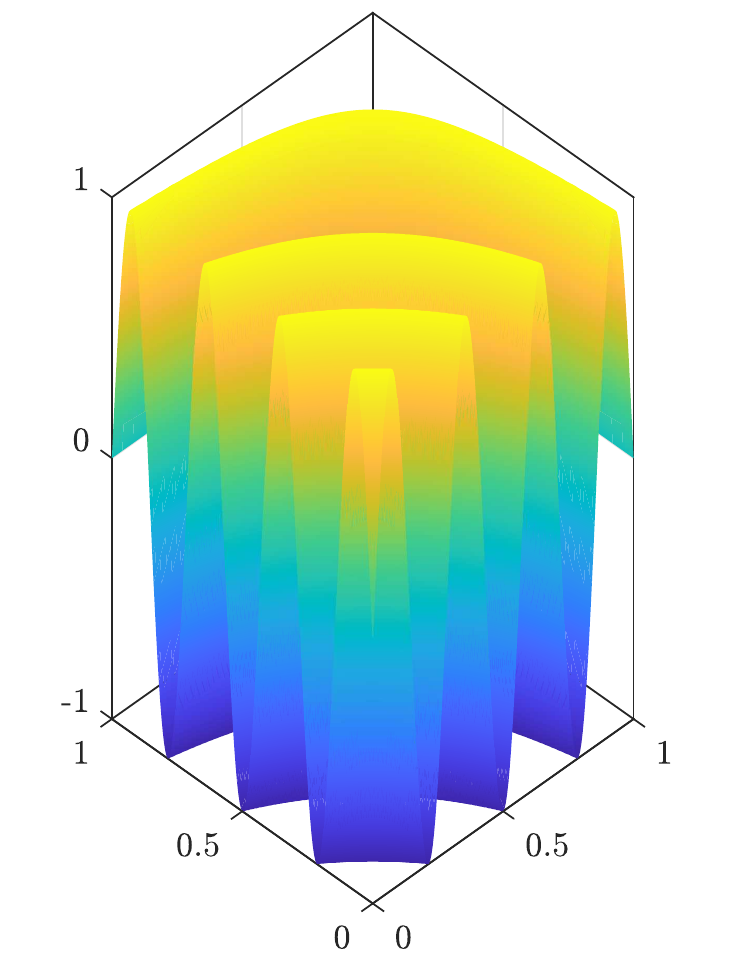}%
\includegraphics[width=0.33\textwidth]{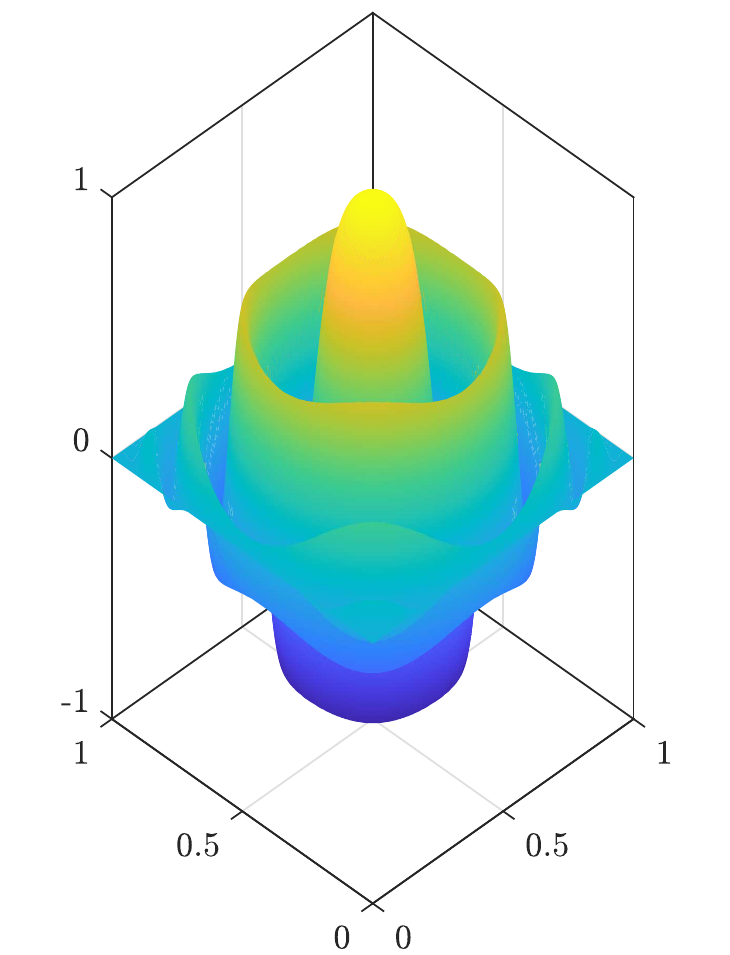}%
\includegraphics[width=0.33\textwidth]{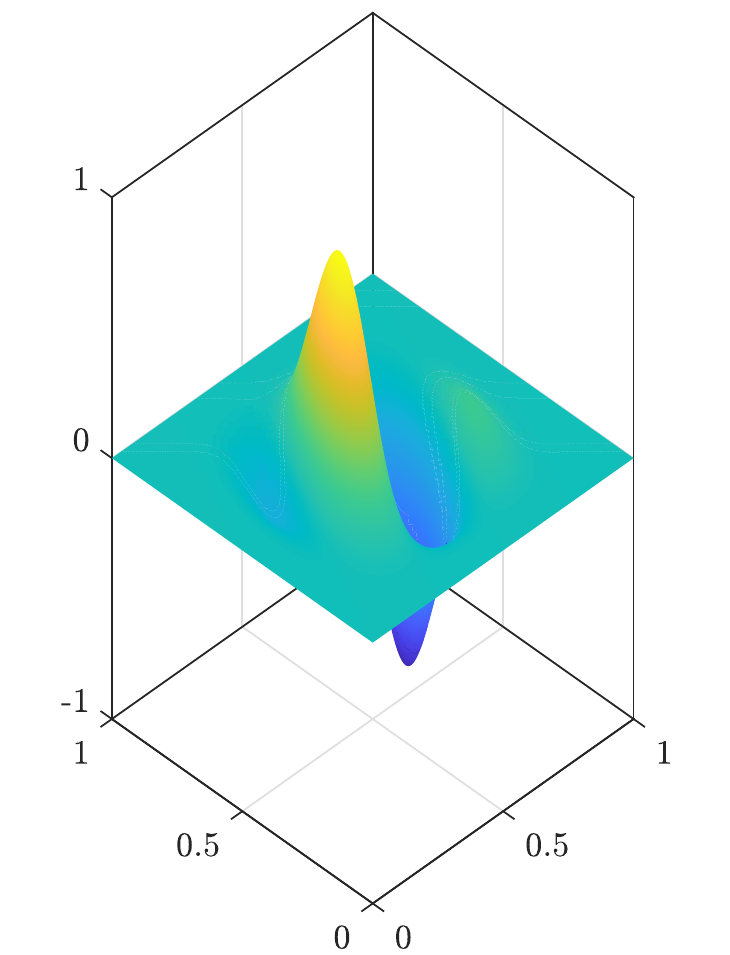}
\caption{Graphs of the functions \eqref{eq:lsqfun}, \eqref{eq:poissonfun}, \eqref{eq:biharmonicfun} (from left to right) used in the numerical examples.}
\label{fig:funs}
\end{figure}

\begin{figure}[t!]
\centering
\includegraphics[width=0.25\textwidth]{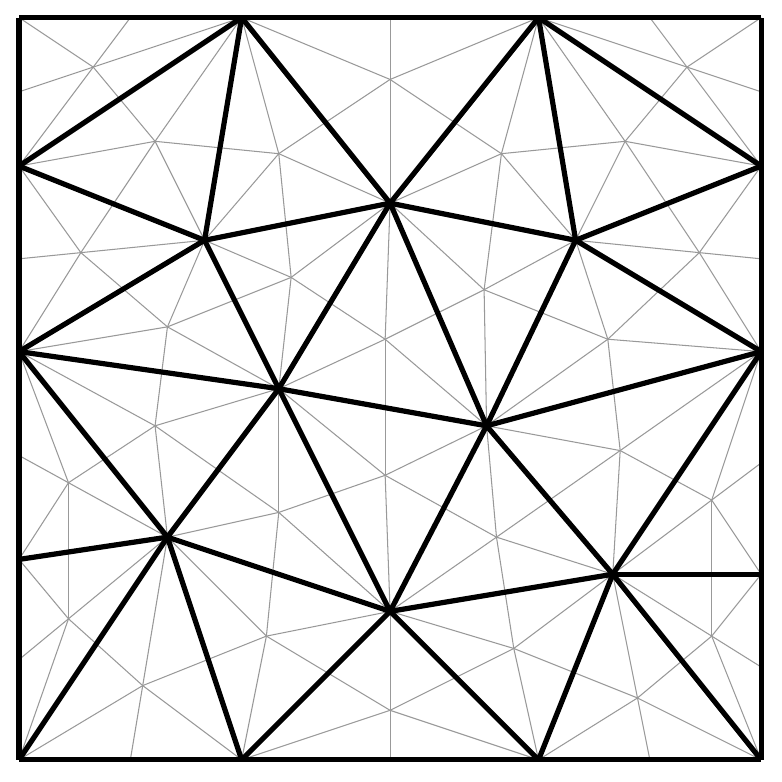}%
\includegraphics[width=0.25\textwidth]{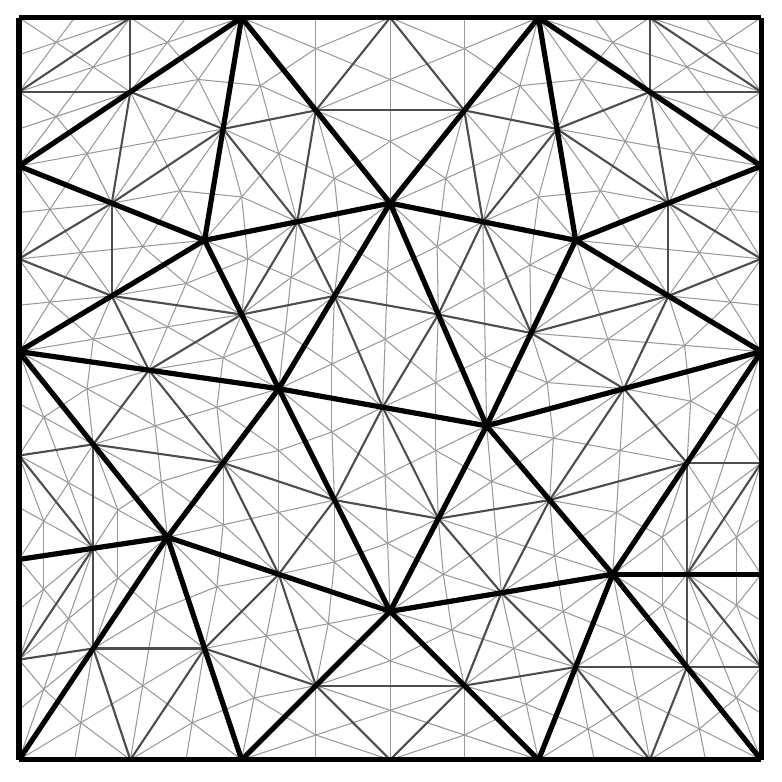}%
\includegraphics[width=0.25\textwidth]{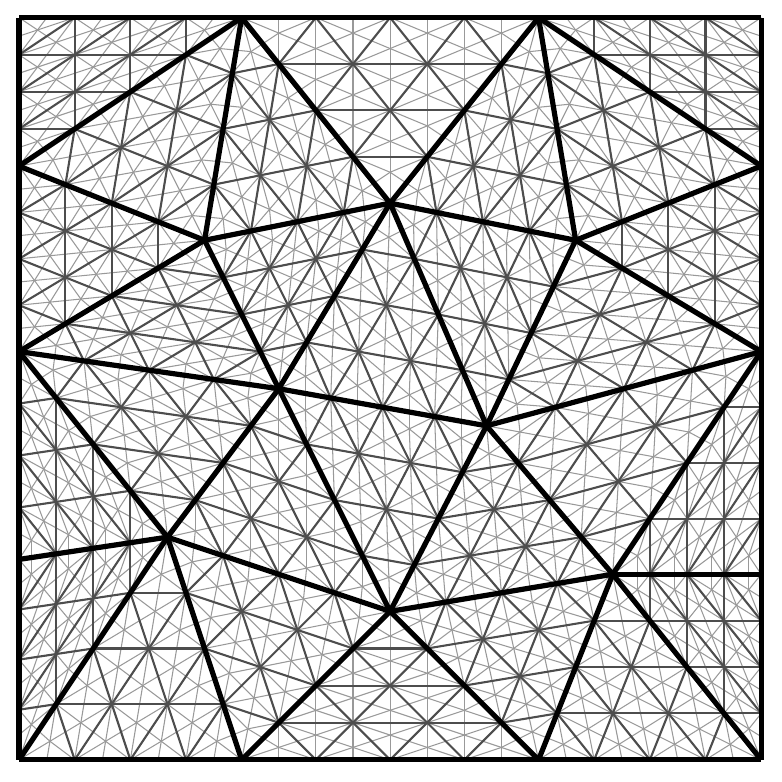}%
\includegraphics[width=0.25\textwidth]{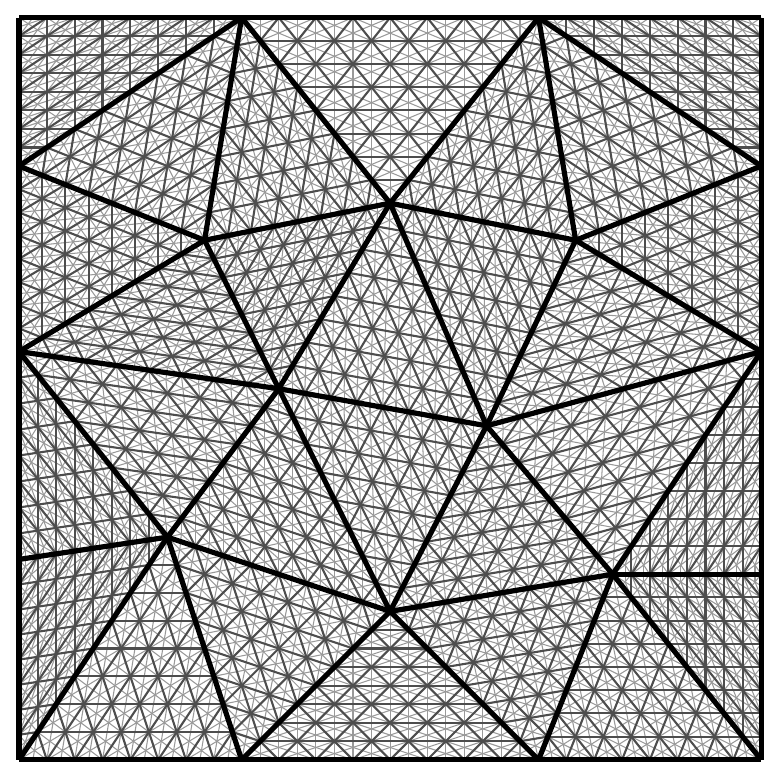}
\caption{An unstructured triangulation $\tri = \tri^1$ of the unit square (left) and three subsequent uniform refinements $\tri^\ell$ for $\ell = 2^k$, $k = 1, 2, 3$ (from left to right). The Powell--Sabin refinements of $\tri^\ell$ are depicted in light gray.}
\label{fig:triangulations}
\end{figure}

\begin{figure}[t!]
\centering
\includegraphics[width=0.33\textwidth]{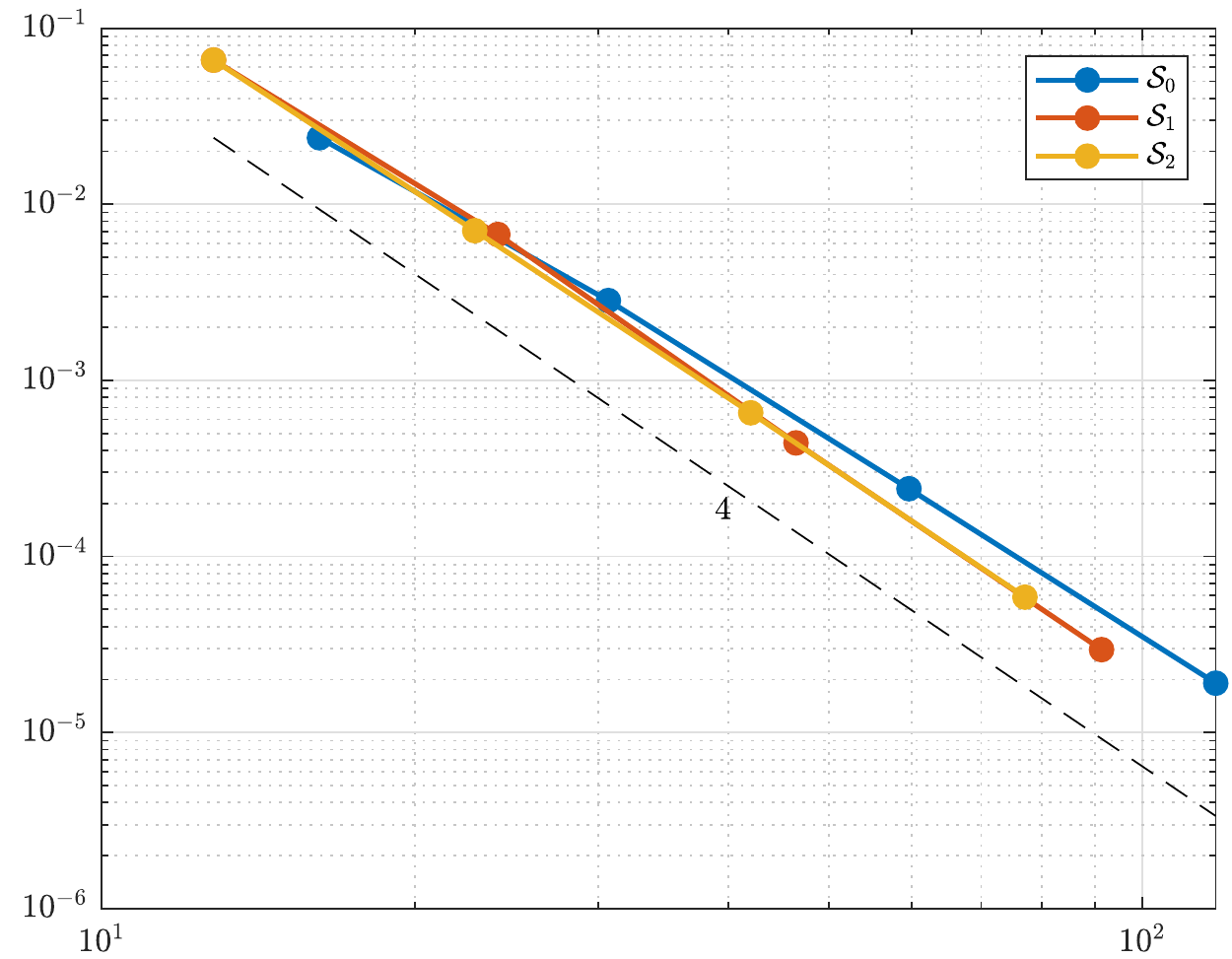}%
\includegraphics[width=0.33\textwidth]{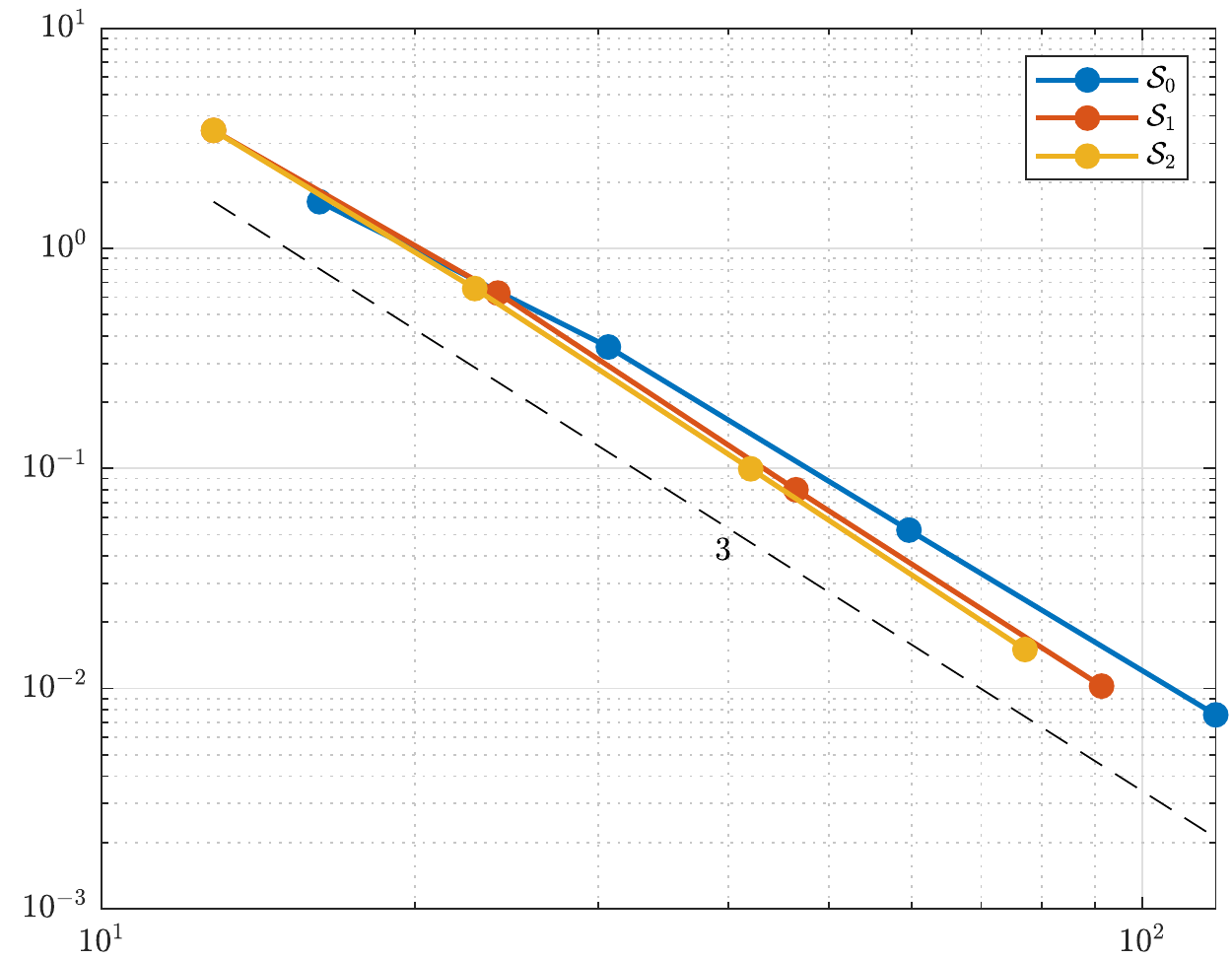}%
\includegraphics[width=0.33\textwidth]{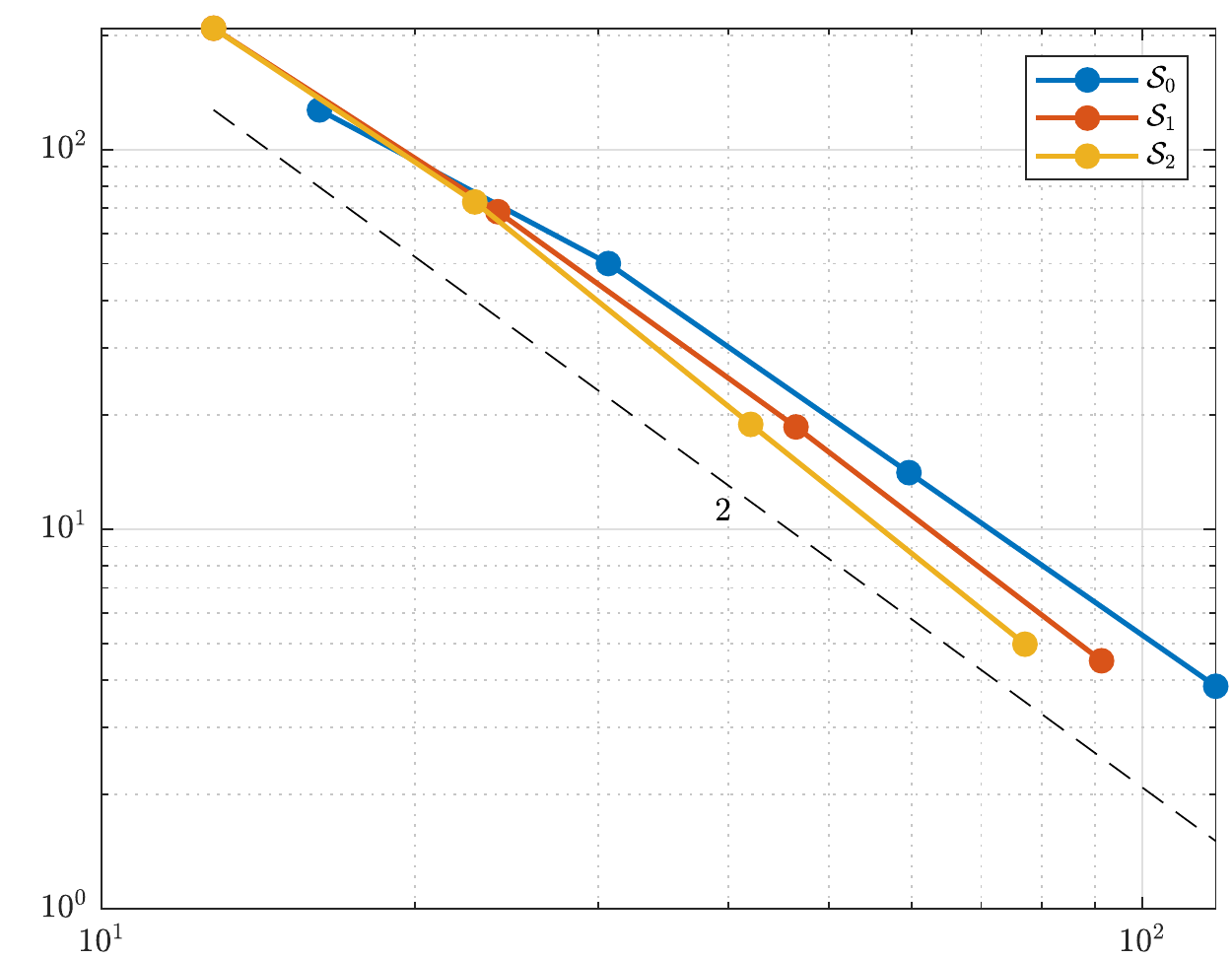}
\caption{$L^2$, $H^1$, and $H^2$ errors (from left to right) with respect to the square root of NDOF in the least squares approximation problem of Example~\ref{ex:lsq}.}
\label{fig:lsq}
\end{figure}

\begin{example}
\label{ex:lsq}
We approximate the function
\begin{equation}
\label{eq:lsqfun}
f(x,y) = \sin\bracket{7 \pi (1-x) (1-y)}
\end{equation}
on $\Omega = [0, 1] \times [0,1]$; see Figure~\ref{fig:funs} (left). The domain is initially partitioned by the triangulation $\tri$ shown in Figure~\ref{fig:triangulations} and subsequently refined by the triangulations $\tri^\ell$ for $\ell = 2^k$, $k = 1, 2, 3$. Figure~\ref{fig:lsq} shows the $L^2$, $H^1$, and $H ^2$ errors of the least squares approximations from the spaces $\splspace_r(\PST^\ell)$, $r = 0, 1, 2$, with respect to the square root of the numbers of degrees of freedom (NDOF) associated with the spaces. As the square root of NDOF is inversely proportional to the length of the longest edge in the triangulation, we observe that the convergence rates are optimal in all three cases with respect to all three errors. We also see that the solutions from the reduced spaces outperform the solutions from the full space in terms of accuracy with respect to NDOF.
\end{example}

\subsection{Approximate solution to a second order boundary value problem}
\label{sec:poisson}

For a function $f \in L^2(\Omega)$, an approximation $s_r$ from $\splspace_r(\PST^\ell)$ to the solution $u$ of the second order boundary value problem
\begin{equation}
\label{eq:poisson}
\left\{
\begin{aligned}
-\nabla^2 
u(p) &= f(p) & \text{if $p$} &\in \Omega, \\
u(p) & = 0 & \text{if $p$} &\in \partial \Omega,
\end{aligned}
\right.
\end{equation}
is commonly obtained by considering the problem in the variational formulation
\begin{equation}
\label{eq:poisson_weak}
\int_\Omega \nabla u(p) \cdot \nabla v(p) \, \mathrm{d}p = \int_\Omega f(p) v(p) \, \mathrm{d}p, \quad v \in H_0^1(\Omega).
\end{equation}
Suppose $\mathring{\map{S}}_r$ is the column vector consisting of the basis functions in $\map{S}_r$ that are zero on the boundary of $\Omega$. When discretizing $H_0^1(\Omega)$ with the space spanned by the components of $\mathring{\map{S}}_r$, the conditions in \eqref{eq:poisson_weak} can be expressed as the system of linear equations
\begin{equation*}
\bracket{ \int_\Omega \bracket{\mathring{\map{M}}_r^{(1,1,1)}(p) + \mathring{\map{M}}_r^{(0,0,1)}(p)} \, \mathrm{d}p } \cdot \mathring{c}_r =
\int_\Omega \mathring{\map{S}}_r(p) \cdot f(p) \, \mathrm{d}p,
\end{equation*}
where $\mathring{\map{M}}_r^{(1,1,1)}$ and $\mathring{\map{M}}_r^{(0,0,1)}$ are defined in the same way as  $\map{M}_r^{(1,1,1)}$ and $\map{M}_r^{(0,0,1)}$ in \eqref{eq:splineprodmatrix} but with $\map{S}_r$ replaced by $\mathring{\map{S}}_r$. The solution $\mathring{c}_r$ of the system determines the approximation $s_r = \mathring{\map{S}}_r^T \cdot \mathring{c}_r$.

\begin{figure}[t!]
\centering
\includegraphics[width=0.33\textwidth]{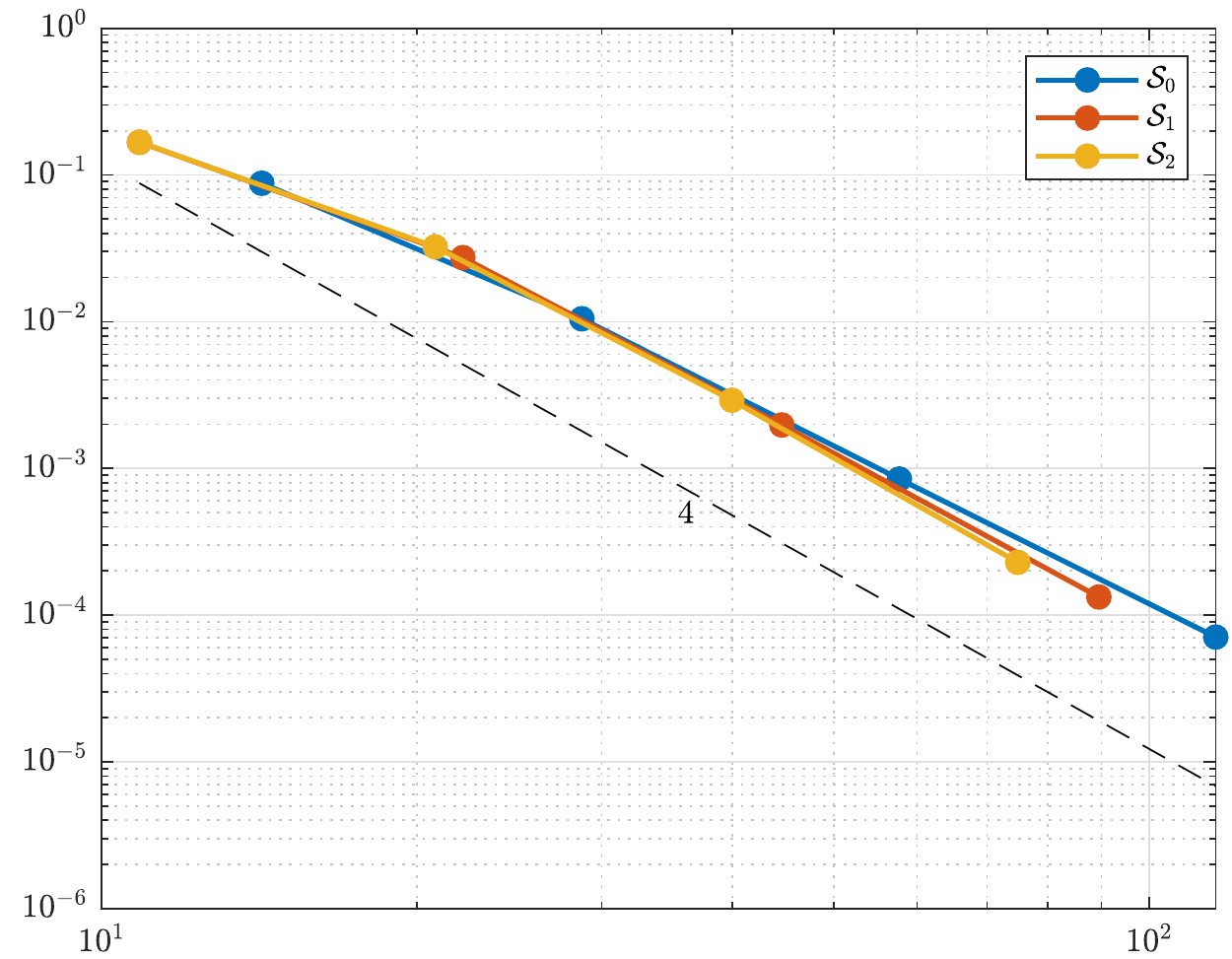}%
\includegraphics[width=0.33\textwidth]{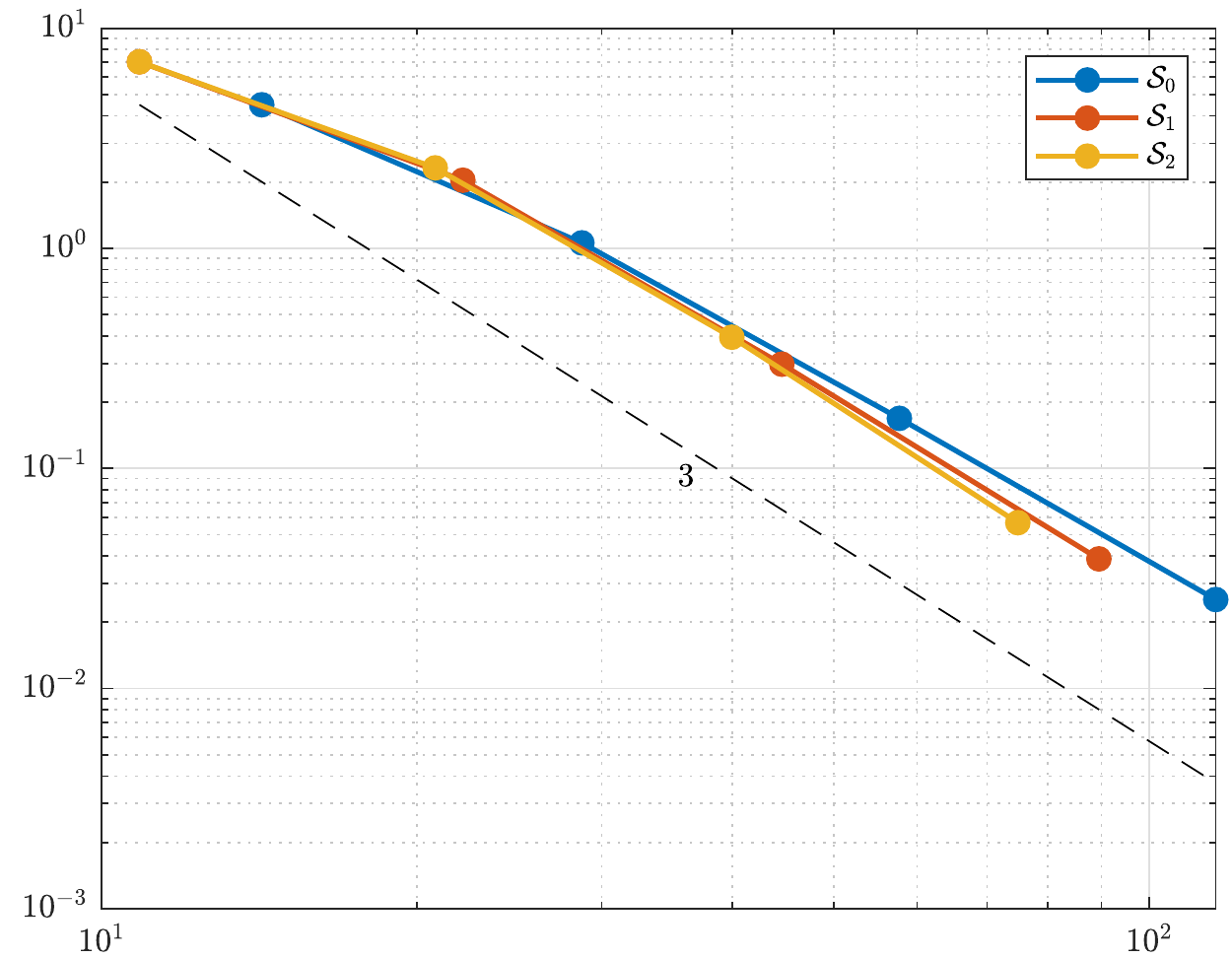}%
\includegraphics[width=0.33\textwidth]{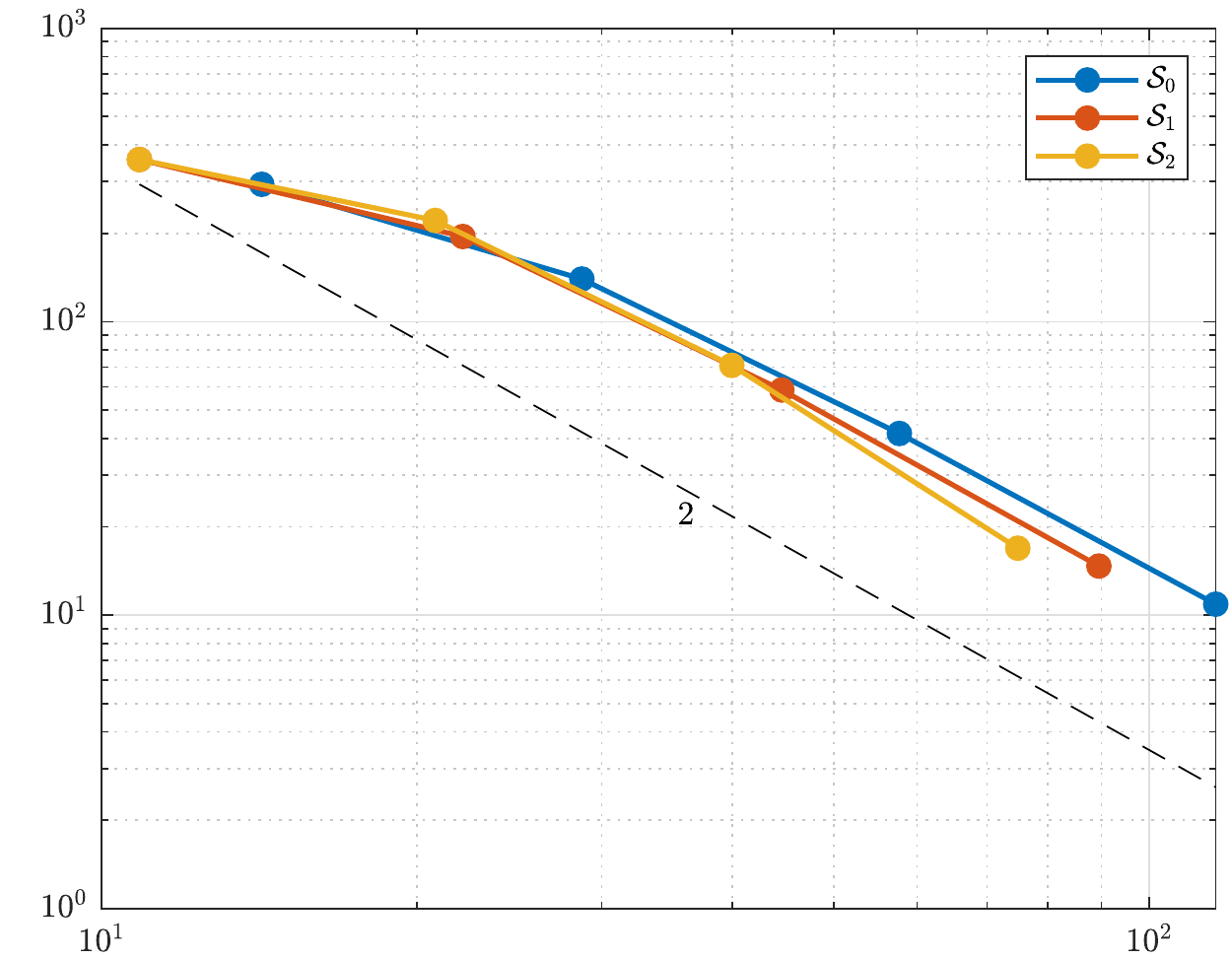}
\caption{$L^2$, $H^1$, and $H^2$ errors (from left to right) with respect to the square root of NDOF in the second order boundary value problem of Example~\ref{ex:poisson}.}
\label{fig:poisson}
\end{figure}

\begin{example}
\label{ex:poisson}
We approximate the manufactured solution
\begin{equation}
\label{eq:poissonfun}
u(x,y) = 16 x (1-x) y (1-y) \cos\bracket{16 \pi \bracket{(x - \tfrac{1}{2})^2+(y - \tfrac{1}{2})^2}},
\end{equation}
shown in Figure~\ref{fig:funs} (middle), to the problem \eqref{eq:poisson} on the square domain partitioned as in Example~\ref{ex:lsq} (Figure~\ref{fig:triangulations}). The $L^2$, $H^1$, and $H^2$ errors of the spline approximations from $\splspace_r(\PST^\ell)$ are provided in Figure~\ref{fig:poisson}. We see again that the convergence rates are optimal and that the solutions from the reduced spaces give better accuracy with respect to NDOF.
\end{example}

\subsection{Approximate solution to a fourth order boundary value problem}

Finally, for $f \in L^2(\Omega)$, we consider the fourth order boundary value problem
\begin{equation}
\label{eq:biharmonic}
\left\{
\begin{aligned}
\nabla^4 
u(p) &= f(p) & \text{if $p$} &\in \Omega, \\
u(p) & = 0 & \text{if $p$} &\in \partial \Omega, \\
\nabla u(p)\cdot\nu &= 0 & \text{if $p$} &\in \partial \Omega,
\end{aligned}
\right.
\end{equation}
where $\nu$ is the unit outer normal to $\partial \Omega$.
The solution $u$ can be described through the variational formulation
\begin{equation*}
\int_\Omega 
\nabla^2 
u(p) \,
\nabla^2 
v(p) \, \mathrm{d}p = \int_\Omega f(p) v(p) \, \mathrm{d}p, \quad v \in H_0^2(\Omega).
\end{equation*}
Let $\mathring{\map{S}}_r$ be now the column vector consisting of the basis functions in $\map{S}_r$ that are zero on the boundary of $\Omega$ and whose normal derivatives are zero on the boundary of $\Omega$. After discretizing $H_0^2(\Omega)$ with the space spanned by the components of $\mathring{\map{S}}_r$, we arrive at the system of linear equations
\begin{equation*}
\bracket{ \int_\Omega \bracket{\mathring{\map{M}}_r^{(2,2,2)}(p) + \mathring{\map{M}}_r^{(2,0,2)}(p) + \mathring{\map{M}}_r^{(0,2,2)}(p) + \mathring{\map{M}}_r^{(0,0,2)}(p)} \, \mathrm{d}p } \cdot \mathring{c}_r =
\int_\Omega \mathring{\map{S}}_r(p) \cdot f(p) \, \mathrm{d}p,
\end{equation*}
where the mappings appearing in the integral on the left side are defined as their analogous counterparts in \eqref{eq:splineprodmatrix} with $\map{S}_r$ replaced by $\mathring{\map{S}}_r$. The solution $\mathring{c}_r$ of the system specifies an approximation $s_r = \mathring{\map{S}}_r^T \cdot \mathring{c}_r$ to $u$ from $\splspace_r(\PST^\ell)$.

\begin{figure}[t!]
\centering
\includegraphics[width=0.33\textwidth]{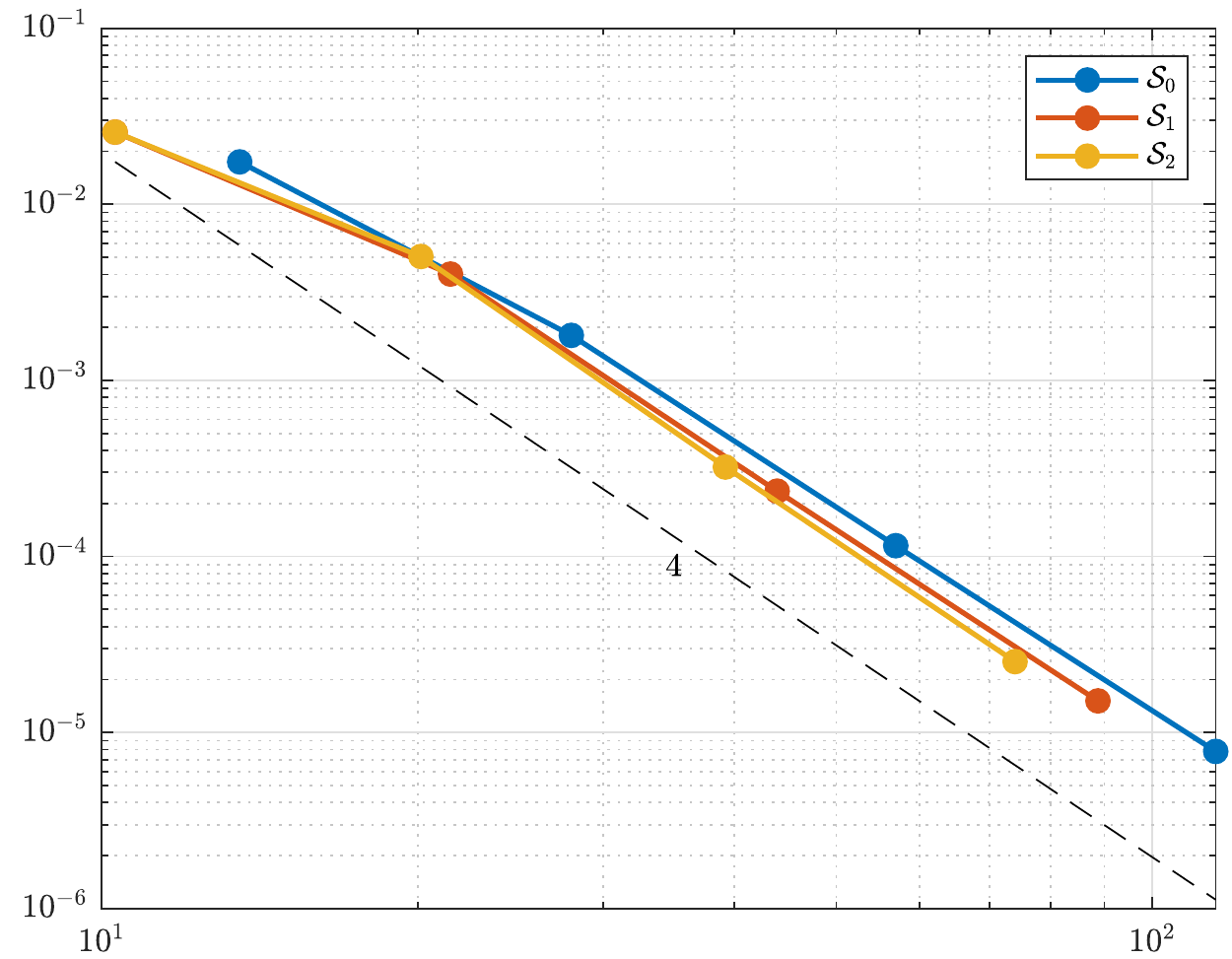}%
\includegraphics[width=0.33\textwidth]{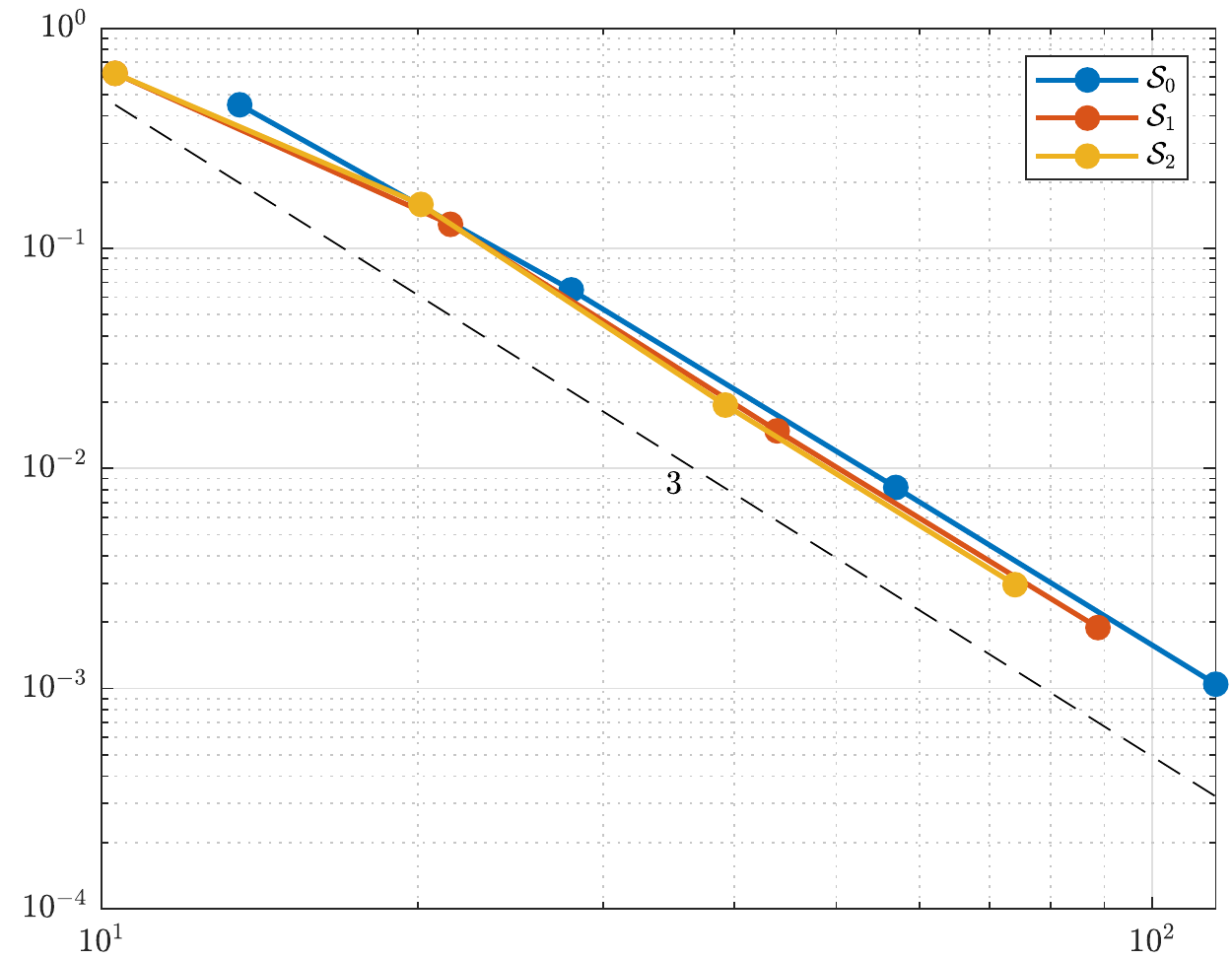}%
\includegraphics[width=0.33\textwidth]{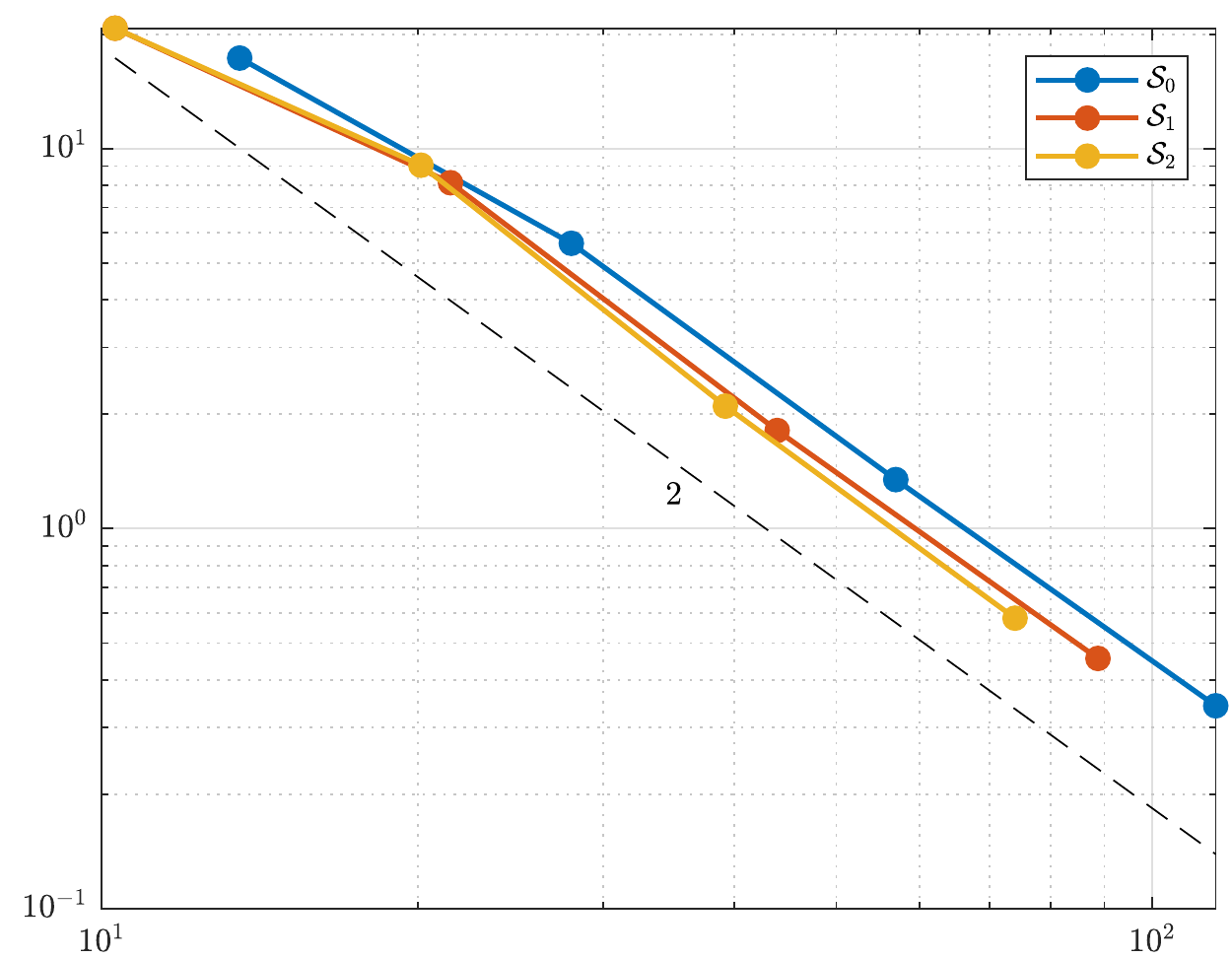}
\caption{$L^2$, $H^1$, and $H^2$ errors (from left to right) with respect to the square root of NDOF in the fourth order boundary value problem of Example~\ref{ex:biharmonic}.}
\label{fig:biharmonic}
\end{figure}

\begin{example}
\label{ex:biharmonic}
As in Examples \ref{ex:lsq} and \ref{ex:poisson}, we consider the square domain partitioned by a sequence of triangulations shown in Figure~\ref{fig:triangulations}. Figure~\ref{fig:biharmonic} shows the $L^2$, $H^1$, and $H^2$ errors of the approximate solutions to the problem \eqref{eq:biharmonic} with the manufactured solution
\begin{equation}
\label{eq:biharmonicfun}
u(x,y) = \sin\bracket{2 \pi (2x-y)} \sin(\pi x)^4 \sin(\pi y)^4
\end{equation}
plotted in Figure~\ref{fig:funs} (right). The results indicate the same behavior of the approximations from the spaces $\splspace_r(\PST^\ell)$ as observed in the previous examples. Indeed, they all exhibit optimal convergence rates, and the reductions prove to be beneficial in achieving better accuracy with respect to NDOF.
\end{example}

\begin{remark}
The basis functions in $\mathring{\map{S}}_r$ used for the fourth order problem here are different from the ones used for the second order problem in Section~\ref{sec:poisson}. Let us elaborate this in more detail for $r=2$. For each boundary edge $e_j \in E^\ell$ and the corresponding triangle $t_k \in T^\ell$ that contains $e_j$, the basis function $B_{\xi_2(e_j, t_k)}^2$ has zero values but not zero normal derivatives on the boundary. An example of such a basis function is depicted in red in Figure~\ref{fig:bsplines12} (right). Moreover, for each boundary vertex $v_i \in V^\ell$ connecting two collinear boundary edges and a proper choice of the triangle $q_i^v$ (with one side on the boundary), one of the three corresponding basis functions also has zero values but not zero normal derivatives on the boundary. An example of such a basis function is depicted in blue in Figure~\ref{fig:bsplines_0} (left).
\end{remark}

\section{Conclusion}
\label{sec:conclusion}

It is very common that the domain over which spline methods are used is in large part partitioned regularly, i.e., only small regions are subjected to irregular meshing. Even though often negligible in area, the irregular parts are far from insignificant in spline constructions and may be, when not properly handled, the cause of decreased accuracy and shape imperfections. In this paper, we have chosen to address such a problem by applying an unstructured spline technology and by reducing the number of degrees of freedom associated with regularly refined regions through the imposition of extra smoothness conditions.

The main contribution is an extraction process based on cubic Powell--Sabin B-splines. Thanks to the linear independence, local support, $C^1$ continuity, non-negativity, partition of unity, and cubic precision, these spline functions provide a sound and reliable approximation framework over any given triangulation. The extraction of smoother functions from the initial set of basis functions does not compromise the optimal approximation power, which in numerical examples manifests in better accuracy with respect to the number of degrees of freedom.

The extraction is presented as a two-step process. Firstly, a partial $C^2$ smoothness is gained inside triangles, i.e., over certain edges and split points of the Powell--Sabin refinement. Secondly, a complete $C^2$ smoothness is obtained over the interior of each triangle with symmetric configuration. The transition between the full set of basis functions and the two reduced alternatives is simple and is expressed in terms of sparse extraction matrices. These can be directly applied in approximation methods, which is in the present paper demonstrated in the context of the least squares approximation method and the finite element method for solving second and fourth order boundary value problems.

Although the paper considers the extraction process over an initial unstructured triangulation that is subsequently uniformly refined, the presented techniques are applicable whenever one can identify the set of triangles with symmetric configuration. Moreover, it is likely that the notion of symmetric configuration can be loosened to a more general geometric symmetry when the edges of the Powell--Sabin refinement meet with only three different slopes at the triangle split point. This extension, which requires a different recombination of basis functions in order to preserve cubic precision, is a subject for future research.

\section*{Acknowledgements}
J.~Gro\v{s}elj was partially supported by the research programme P1-0294 of Javna agencija za raziskovalno dejavnost Republike Slovenije (ARRS). H.~Speleers was supported in part by a GNCS 2022 project (CUP E55F22000270001) of Gruppo Nazionale per il Calcolo Scientifico -- Istituto Nazionale di Alta Matematica (GNCS -- INdAM) and by the MUR Excellence Department Project MatMod@TOV (CUP E83C23000330006) awarded to the Department of Mathematics of the University of Rome Tor Vergata.

\appendix
\section{Cubic polynomial reproduction}
\label{sec:polynomials}

In this appendix, we show that the second reduced spline space $\splspace_2(\PST^\ell)$ contains the space $\polspace_3$ of bivariate polynomials of total degree less than or equal to $3$.
The characterization of the first reduced spline space $\splspace_1(\PST^\ell)$ given by \cite{ps3_speleers_15} tells us that any $g \in \polspace_3$ belongs to $\splspace_1(\PST^\ell)$. Considering \eqref{eq:spline},
it can be uniquely represented as
\begin{equation*}
g(p) =  \map{S}_1(p)^T \cdot c_1
\end{equation*}
for some column vector $c_1 \in \R^{m+m_1}$.
Following the reduction procedure from Section~\ref{sec:reduction2}, we see that the basis functions $B_1, B_2, \ldots, B_m$ remain unchanged in the basis $\map{S}_2$. The same holds true for the basis functions $B_{\xi_1(e_j,t_k)}^1$, $B_{\xi_1(e_j, t_{k'})}^1$ for each pair $(e_j, t_k), (e_j, t_{k'}) \in \Xi_1$ such that both $t_k$ and $t_{k'}$ are in $T^\ell \setminus \Tsym^\ell \cup \{ \emptyset \}$. Therefore, in order to verify that $g(p)$ can be represented in terms of $\map{S}_2(p)$, we only need to focus on the remaining set of basis functions.

Let $t_k \in \Tsym^\ell$ be a triangle with vertices $v_{i_1}$, $v_{i_2}$, and $v_{i_3}$. Moreover, let $t_{k'} \in T^\ell$ be the neighboring triangle with vertices $v_{i_1}$, $v_{i_2}$, and $v_{i_3'}$. Finally, let $e_j \in E^\ell$ be the edge between the vertices $v_{i_1}$ and $v_{i_2}$. Since $t_k$ has symmetric configuration, we know that
\begin{align*}
v_k^t &= \frac{v_{i_1} + v_{i_2} + v_{i_3}}{3} = \frac{2v_{i_3} + v_{i_3'}}{3}, \\
v_{k'}^t &= \frac{v_{i_1} + v_{i_2} + v_{i_3'}}{3} = \frac{v_{i_3} + 2v_{i_3'}}{3}.
\end{align*}
Let $\bloss[g]$ denote the blossom of the polynomial $g \in \polspace_3$; see \cite{seidel_93,bb_speleers_11} for more details about blossoming. 
From \cite{ps3_groselj_17} we know that $g(p)$ can be represented in terms of $\map{S}_1(p)$ as
\begin{equation*}
g(p) = \bloss[g](v_{i_1}, v_{i_2}, v_k^t) B_{\xi_1(e_j, t_k)}^1(p) + \bloss[g](v_{i_1}, v_{i_2}, v_{k'}^t) B_{\xi_1(e_j, t_{k'})}^1(p) + R_1(p),
\end{equation*}
where $R_1(p)$ collects the contributions of the remaining basis functions in $\map{S}_1(p)$.
By exploiting the multi-affinity property of blossoms, we obtain
\begin{align*}
g(p) &= \bloss[g]\Bigl(v_{i_1}, v_{i_2}, \frac{2v_{i_3} + v_{i_3'}}{3}\Bigr) B_{\xi_1(e_j, t_k)}^1(p) + \bloss[g]\Bigl(v_{i_1}, v_{i_2}, \frac{v_{i_3} + 2v_{i_3'}}{3}\Bigr) B_{\xi_1(e_j, t_{k'})}^1(p) + R_1(p) \\
&= \bloss[g](v_{i_1}, v_{i_2}, v_{i_3}) \Bigl[ \frac{2}{3}B_{\xi_1(e_j, t_k)}^1(p) + \frac{1}{3}B_{\xi_1(e_j, t_{k'})}^1(p)\Bigr] \\
&\quad + \bloss[g](v_{i_1}, v_{i_2}, v_{i_3'}) \Bigl[ \frac{1}{3}B_{\xi_1(e_j, t_k)}^1(p) + \frac{2}{3}B_{\xi_1(e_j, t_{k'})}^1(p)\Bigr] + R_1(p).
\end{align*}
After repeating the above argument for all pairs of triangles with at least one triangle in $\Tsym^\ell$, we arrive at
\begin{equation*}
g(p) = \bloss[g](v_{i_1}, v_{i_2}, v_{i_3}) B_{\xi_2(\emptyset, t_k)}^2(p) + \bloss[g](v_{i_1}, v_{i_2}, v_{i_3'}) B_{\xi_2(\emptyset, t_{k'})}^2(p) + R_2(p)
\end{equation*}
if $t_{k'} \in \Tsym$, or
\begin{equation*}
g(p) = \bloss[g](v_{i_1}, v_{i_2}, v_{i_3}) B_{\xi_2(\emptyset, t_k)}^2(p) + \bloss[g](v_{i_1}, v_{i_2}, v_{i_3'}) B_{\xi_2(e_j, t_{k'})}^2(p) + R_2(p)
\end{equation*}
if $t_{k'} \in T^\ell \setminus \Tsym^\ell$, where in both cases $R_2(p)$ collects the contributions of the remaining basis functions in $\map{S}_2(p)$. This means that every polynomial $g \in \polspace_3$ can be represented in the basis $\map{S}_2$ and so belongs to the space $\splspace_2(\PST^\ell)$.

\bibliography{references}

\end{document}